\documentclass[reqno,12pt]{amsart}
\usepackage{amsfonts,amsmath,amsthm,amssymb,amscd}
\newcommand{\version}{Ver.~0.0}
\newcommand{\setversion}[1]{\renewcommand{\version}{Ver.~{#1}}}
\setversion{1.0$\alpha$ [Wed Jul 17 10:06:19 JST 2002]}
\setversion{2.0 [Tue Sep  9 17:54:13 JST 2003]}
\setversion{2.1 [Wed Sep 10 16:25:07 JST 2003]}
\setversion{3.0 [Fri Oct  3 18:26:10 JST 2003]}
\setversion{3.1 [Sun Dec  7 00:39:14 JST 2003]}
\setversion{3.2 [Thu Dec 25 11:25:56 JST 2003]}

\title [lifting of nilpotent orbits]
{Theta lifting of nilpotent orbits for symmetric pairs}

\author{Kyo Nishiyama}
\address{
Department of Mathematics\\
Graduate School of Science\\
Kyoto University\\
Sakyo, Kyoto 606-8502, Japan}
\email{kyo@math.kyoto-u.ac.jp}

\author{Hiroyuki Ochiai}
\address{
Department of Mathematics\\
Nagoya University\\
Nagoya, 464-8602, Japan}
\email{ochiai@math.nagoya-u.ac.jp}

\author{Chen-bo Zhu}
\address{Department of Mathematics\\
National University of Singapore\\
2 Science Drive 2\\
Singapore 117543}
\email{matzhucb@nus.edu.sg}


\theoremstyle{plain}
\newtheorem{theorem}{Theorem}
\newtheorem{proposition}[theorem]{Proposition}
\newtheorem{corollary}[theorem]{Corollary}
\newtheorem{lemma}[theorem]{Lemma}

\newtheorem{introtheorem}{Theorem}

\newtheorem{mytheorem}{Theorem}
\newtheorem{myproposition}[mytheorem]{Proposition}

\newtheorem{mylemma}[mytheorem]{Lemma}

\theoremstyle{definition}
\newtheorem{definition}[theorem]{Definition}
\theoremstyle{remark}
\newtheorem{remark}[theorem]{\upshape Remark}

\numberwithin{equation}{section}
\numberwithin{theorem}{section}

\newcommand{\Z}{\mathbb{Z}}
\newcommand{\Q}{\mathbb{Q}}
\newcommand{\R}{\mathbb{R}}

\newcommand{\bbK}{\mathbb{K}}
\newcommand{\bbG}{\mathbb{G}}
\newcommand{\C}{\mathbb{C}}
\newcommand{\bbP}{\mathbb{P}}

\newcommand{\lie}[1]{\mathfrak{#1}}

\newcommand{\liec}[1]{\mathfrak{#1}}

\newcounter{thmenum}
\newenvironment{thmenumerate}{%
\begin{list}{$(\thethmenum)$}{%
\usecounter{thmenum}
\setlength{\labelsep}{.5em}
\setlength{\labelwidth}{0pt}
\setlength{\parsep}{0pt}
\setlength{\leftmargin}{3pt}
\setlength{\rightmargin}{0pt}
\setlength{\itemindent}{\leftmargin}
}}
{\end{list}}

\newcommand{\mycomment}[1]{} 

\newlength{\lengthcup}
\newcommand{\vertin}{\makebox[0pt][l]{$\cup$}\hspace*{.36\lengthcup}\raisebox{.31ex}{\scriptsize $|$}}

\newcommand{\mathoperator}[1]{\mathop\mathrm{#1}\nolimits{}}

\newcommand{\transpose}[1]{\,{}^t{#1}}

\newcommand{\Hom}{\mathop\mathrm{Hom}\nolimits{}}

\newcommand{\rank}{\mathop\mathrm{rank}\nolimits{}}
\newcommand{\codim}{\mathop\mathrm{codim}\nolimits{}}

\newcommand{\directsum}{\mathop{\ \sum\nolimits^{\oplus}}}

\renewcommand{\Im}{\mathop\mathrm{Im}\nolimits{}}

\newcommand{\bbI}{\mathbb{I}}

\newcommand{\closure}[1]{\overline{#1}}

\newcommand{\Omin}{\mathbb{O}_{\mathrm{min}}}
\newcommand{\trivial}{\mathbf{1}}

\newcommand{\irreps}[1]{\mathrm{Irr}(#1)}
\newcommand{\restrict}{\big|}
\newcommand{\partition}{\mathcal{P}}

\newcommand{\Spec}{\mathop\mathrm{Spec}\nolimits{}}
\newcommand{\Sym}{\mathop\mathrm{Sym}\nolimits{}}
\newcommand{\Alt}{\mathop\mathrm{Alt}\nolimits{}}

\newcommand{\orbit}{\mathbb{O}}  %

\newcommand{\nullcone}{\lie{N}}

\newcommand{\nilpotents}{\mathcal{N}}

\newcommand{\theoremref}[1]{Theorem \ref{#1}}

\newcommand{\corollaryref}[1]{Corollary \ref{#1}}

\newcommand{\eqnref}[1]{{\upshape (\ref{#1})}}

\newcommand{\harmonics}{\mathcal{H}}
\newcommand{\harmonicsp}{\harmonics^{+}}
\newcommand{\harmonicsq}{\harmonics^{-}}

\newcommand{\composit}{\odot}

\newcommand{\Thetalift}{\theta}

\newcommand{\bbKC}{\mathbb{K}_{\C}}

\newcommand{\Gbig}{G}
\newcommand{\Kbig}{K}
\newcommand{\KbigC}{K_{\C}}
\newcommand{\Kbigp}{K^{+}}
\newcommand{\Kbigq}{K^{-}}
\newcommand{\Kbigpq}{K^{\pm}}
\newcommand{\KbigpC}{K_{\C}^{+}}
\newcommand{\KbigqC}{K_{\C}^{-}}
\newcommand{\KbigpqC}{K_{\C}^{\pm}}
\newcommand{\Lbig}{L}
\newcommand{\Lbigp}{L^{+}}
\newcommand{\Lbigq}{L^{-}}
\newcommand{\LbigpC}{L_{\C}^{+}}

\newcommand{\Mbig}{M}
\newcommand{\Gsmall}{G'}
\newcommand{\Ksmall}{K'}
\newcommand{\KsmallC}{K'_{\C}}
\newcommand{\Lsmall}{L'}
\newcommand{\LsmallC}{L'_{\C}}
\newcommand{\Msmall}{M'}

\newcommand{\ssmall}{\lie{s'}}
\newcommand{\ssmallp}{\lie{s'_{+}}}
\newcommand{\ssmallq}{\lie{s'_{-}}}
\newcommand{\ssmallpq}{\lie{s'_{\pm}}}
\newcommand{\gsmall}{\liec{g'}}
\newcommand{\ksmall}{\liec{k'}}
\newcommand{\sbig}{\lie{s}}
\newcommand{\gbig}{\liec{g}}
\newcommand{\kbig}{\liec{k}}
\newcommand{\Lie}{\mathop\mathrm{Lie}\nolimits{}}

\newcommand{\momentbig}{\varphi}
\newcommand{\momentsmall}{\psi}

\newcommand{\momentsmallq}{\psi^{-}}
\newcommand{\momentsmallpq}{\psi^{\pm}}

\newcommand{\Wp}{W^{+}}
\newcommand{\Wq}{W^{-}}
\newcommand{\Wpq}{W^{\pm}}
\newcommand{\nullconep}{\nullcone^{+}}
\newcommand{\nullconeq}{\nullcone^{-}}
\newcommand{\nullconepq}{\nullcone^{\pm}}

\newcommand{\Obig}{\orbit}
\newcommand{\Osmall}{{\orbit}'}

\newcommand{\sorbit}{\mathcal{O}}
\newcommand{\sorbitp}{\mathcal{O}^+}
\newcommand{\sorbitq}{\mathcal{O}^-}
\newcommand{\sorbitpq}{\mathcal{O}^{\pm}}

\newcommand{\Osmalltrivial}{\orbit'{}^{\trivial}}
\newcommand{\Obigtrivial}{\orbit{}^{\trivial}}

\newcommand{\Obighol}{\orbit{}^{\hslash}}

\newcommand{\GITquotient}{/\!/}

\begin{document}

\begin{abstract} 
We consider a reductive dual pair $(\Gbig, \Gsmall)$ in the stable range
with $\Gsmall $ the smaller member and of Hermitian symmetric type. 
We study the theta lifting of nilpotent $\KsmallC $-orbits, where $ \Ksmall $ is a maximal compact subgroup of $ \Gsmall $ and we describe the precise $ \KbigC $-module structure of the regular function ring of the closure of the lifted nilpotent orbit of the symmetric pair $ ( \Gbig , \Kbig ) $.  
As an application, we prove sphericality and normality of the closure of certain nilpotent $ \KbigC $-orbits obtained in this way. We also give integral formulas for their degrees. 
\end{abstract}

\maketitle

\tableofcontents
\newpage

\section*{Introduction}

Let $ ( \Gbig , \Gsmall ) $ be a reductive dual pair in a symplectic group $ \bbG = Sp( 2N , \R ) $, where $ N $ denotes the rank of $ \bbG $.  
Throughout this paper, we will assume that $ ( \Gbig , \Gsmall ) $ is of type I, and 
it is in the stable range with $ \Gsmall $ the smaller member (cf. \cite{Howe.1979, Li.1989}). 
We will also assume that 
$ \Gsmall / \Ksmall $ is an irreducible Hermitian symmetric space.  
By the classification of irreducible dual pairs, our restriction amounts to saying that $ ( \Gbig , \Gsmall ) $ is in the following list 
(cf.\  \cite{Howe.1989.transcending}).

\begin{table}[htbp]
\caption{The dual pairs treated in this paper}
\label{table:dual.pairs.treated}
$ \displaystyle
\begin{array}{l@{\qquad}l}
\text{the pair $ ( \Gbig , \Gsmall ) $} & \text{stable range condition} \\
\hline
\\
( O( p, q ) , Sp( 2 n, \R ) ) & 2 n < \min (p, q)  \\[5pt]
( U( p, q ), U( m, n ) ) ) & m + n \leq \min (p, q)  \\[5pt]
( Sp( p, q ) , O^{\ast}( 2 n ) ) & n \leq \min (p, q)  
\end{array}
$
\end{table}
\noindent
Note that we have excluded the equality $2n=\min(p,q)$ in the first case, to avoid some small technicalities.
\mycomment{
This condition is equivalent to the condition that all the constituents are holomorphic discrete series 
for the compact dual pairs $ \Kbigpq \times \Gsmall $.  
}

Let $ \gbig $ (resp.\  $ \gsmall $) be the complexification of the Lie algebra of $ G $ (resp.\  $ G' $) 
and $ K $ (resp.\  $ K' $) a maximal compact subgroup of $ G $ (resp.\  $ G' $).  
Let  
\begin{displaymath}
\gbig   = \kbig   \oplus \sbig   , \qquad 
\gsmall = \ksmall \oplus \ssmall , 
\end{displaymath}
be respectively Cartan decompositions of $ \gbig $ and $ \gsmall $, compatible in certain way 
(see \S \ref{sec:lifting.nilpotent.orbit}). 

\medskip
\noindent
\textbf{Notation.}  
Let $ G $ be a reductive algebraic group acting on an affine variety $ X $.  
Then we denote the ring of $ G $-invariants on the coordinate ring by $ \C[ X ]^G $,  
and the affine quotient of $ X $ by $ G $ is defined to be 
\begin{math}
X \GITquotient G = \Spec \C[ X ]^G , 
\end{math}
which parameterizes closed orbits in $ X $.  
The affine variety $ X \GITquotient G $ is often called a categorical quotient.

\medskip
Let $ W_{\R} \simeq \R^{2 N} $ be a real symplectic space which realizes $ \bbG = Sp(2N,\R ) $ as a symplectic group on $ W_{\R} $.  
There is a canonical complex structure on $ W_{\R} $ and we can view $ W_{\R} $ as (the underlying real vector space of) a complex vector space $ W \simeq \C^N $.  
The symplectic form on $ W_{\R} $ is then given by the imaginary part of a canonical positive definite Hermitian form on $ W $.  
By this identification, a maximal compact subgroup $ \bbK $ of $ \bbG $ is realized as the unitary group $ U( W ) \simeq U( N ) $ on $ W \simeq \C^N $.

We may choose maximal compact subgroups $ \Kbig $ and $ \Ksmall $ of $ \Gbig $ and $ \Gsmall $ respectively, in such a way that $ \Kbig \cdot \Ksmall $ is contained in the standard maximal compact subgroup $ \bbK \simeq U(N) $ of $ \bbG $.  

There is a double fibration of $ W $ by the moment maps $ \momentbig $ and $ \momentsmall $ :
\begin{equation*}
\sbig \xleftarrow{\quad \momentbig \quad} W \xrightarrow{\quad \momentsmall \quad} \ssmall .
\end{equation*}
For explicit realization of these maps, see the Appendix.  
Here we only mention some of the important properties of $ \momentbig $ and $ \momentsmall $.  
Thanks to the stable range condition, $ \momentsmall $ is a surjective, affine quotient map by the action of $ \KbigC $.  
Thus we have $ \ssmall \simeq W \GITquotient \KbigC $. Furthermore $ \momentsmall : W \rightarrow \ssmall $ is flat. 
On the other hand, the image of $ \momentbig $ is an irreducible closed subvariety of $ \sbig $ and 
$ \momentbig : W \rightarrow \momentbig( W ) \subset \sbig $ is also an affine quotient map by $ \KsmallC $.  

\begin{introtheorem}[Theorem \ref{thm:fiber.of.Osmall}]
Take a nilpotent $ \KsmallC $-orbit $ \Osmall \subset \nilpotents( \ssmall ) $.  
Then the scheme theoretic fiber 
$ \momentsmall^{-1}( \closure{\Osmall} ) = W \times_{ \ssmall } \closure{ \Osmall } $ 
is a reduced, closed irreducible affine subvariety of $ W $.  
Moreover, it is the closure of a single $ \KbigC \times \KsmallC $-orbit in $ W $.  
\end{introtheorem}

The above result easily implies that the set $ \momentbig( \momentsmall^{-1} ( \closure{\Osmall} ) ) $ is the closure of a single nilpotent $ \KbigC $-oribt $ \Obig $ in $ \sbig $.  
Thus we have a mapping from the set of nilpotent $ \KsmallC $-orbits in $ \ssmall $ to the set of nilpotent $ \KbigC $-orbits in $ \sbig $, which sends $ \Osmall $ to $ \Obig $.  
We call this mapping the theta lifting and denote it by $ \Obig = \theta( \Osmall ) $. In the language of the signed Young diagrams, $ \Obig $ is obtained from $ \Osmall $ by adding one extra box to each row at the right end (see Proposition \ref{prop:signed.Y.diagram}).  Thus $ \Obig $ is slightly ``bigger'' than $ \Osmall $.

We remark that T.~Ohta has proved the irreducibility and the existence of an open dense orbit by a totally different (but case-by-case) method (\cite{Ohta.preprint1}, \cite{Ohta.preprint2}). Another recent work by Daszkiewicz, Kra\'{s}kiewicz and Przebinda \cite{DKP.preprint} proves similar results by classifying $ \KbigC \times \KsmallC $-orbit in $ W $ combinatorially. There have also been earlier investigations of correspondence of complex nilpotent orbits in \cite{DKP.complex} and by S.-Y. Pan for complex as well as real nilpotent orbits. 

Our approach is in some sense intrinsic and it is based on methods of algebraic geometry. More importantly, we hope that it will help to shed light on the geometric aspects of the theta correspondence. The recent work of two of the authors \cite{Nishiyama.Zhu.2001,Nishiyama.Zhu.2003} are efforts in this direction. 

Let $\C [\closure{\Obig}]$ (resp.\  $\C [\closure{\Osmall}]$) be the regular function ring
on the closure of $\Obig$ (resp.\  $\Osmall $). 
We are interested in the $ \KbigC $-module structure of $\C [\closure{\Obig}]$. One reason for this interest comes from Vogan's philosophy of unipotent representations \cite{Vogan.1991}. Very roughly speaking, it entails that if a unitary representation $\pi $ is to be ``associated" to a nilpotent orbit $\Obig$, the regular function ring $\C [\closure{\Obig}]$ should encode the $K$-structure of the representation $\pi $. 
In the present situation, we will show that the $\KbigC$-module structure of $\C [\closure{\Obig}]$ can be described entirely and explicitly by the regular function ring $\C [\closure{\Osmall}]$ and the space of $\KbigC$-harmonic polynomials $ \harmonics(\KbigC) $ on $W$.  
The main results are summarized in Theorem \ref{thm:fun.ring.of.lifted.orbit}.  
In particular we have the following 

\begin{introtheorem}
We have a $ \KbigC $-module isomorphism 
\begin{displaymath}
 \C[ \closure{\Obig} ] \simeq ( \harmonics(\KbigC) \otimes \C[ \closure{\Osmall} ] )^{\KsmallC} .
\end{displaymath}
\end{introtheorem}

It is well-known \cite{Howe.1989.classical} that $\harmonics( \KbigC )$ is multiplicity-free as a representation of
$\KbigC \times \KsmallC$ and the decomposition of $\harmonics( \KbigC )\restrict _{\KbigC \times \KsmallC}$
determines a one-to-one correspondence between irreducible representations of $\KbigC $ and $\KsmallC$ occurring in
$\harmonics( \KbigC )$. We note that the form of isomorphism $\C[ \closure{\Obig} ] \simeq ( \harmonics(\KbigC) \otimes \C[ \closure{\Osmall} ] )^{\KsmallC}$ bears strong resemblance to theta correspondence of admissible representations of $\Gsmall$ and $\Gbig$ \cite{Howe.1989.transcending}. 

In \S \ref{section:lift.hol.nilpotent.orbit}, we will give some examples of the lifted orbits and their regular function rings.  We obtain a family of spherical nilpotent orbits 
and a precise $ \KbigC $-module structure of $ \C[ \closure{\Obig} ] $.  
We note that they occupy a large portion of the set of all the spherical nilpotent $ \KbigC $-orbits, 
which are classified by D.~R.~King \cite{King.preprint} recently.  
As an application we prove normality of the closure of some of these nilpotent orbits, which are lifted from the smaller group.  We also give explicit integral formulas of their degrees. The result is summarized in Proposition \ref{prop:degree.of.orbits.integral.expression}.  
The formulas involve certain integral expression, which can be explicitly evaluated.   

\begin{introtheorem}
Let $\kappa$ be a complex number with a positive real part.
Then
\begin{equation*}
\frac{1}{n!} \int_{\Omega_n} D_n(x^2) D_n(x) \prod_{i=1}^n x_i^{\kappa-1} dx
= \frac{2^{n(n-1)/2} \prod_{i=0}^{n-1} i! \, \Gamma(\kappa+2i)}%
{\Gamma(3n(n-1)/2 + n\kappa +1)} ,
\end{equation*}
where $ D_n(x) = \prod_{1 \leq i < j \leq n} ( x_i - x_j ) $, $ D_n(x^2) = \prod_{1 \leq i < j \leq n} ( x^2_i - x^2_j ) $, 
and the region of integral is given as $ \Omega_n = \{ x = ( x_i )_{1 \leq i \leq n} \mid x_i \geq 0 , \sum_{i = 1}^n x_i \leq 1 \} $.
\end{introtheorem}

In the Appendix, we describe the explicit form of moment maps and give some basic properties of 
null cones.  The results in the Appendix are known to the experts, but we include them because of lack of appropriate references.

\section{Diamond pairs}
\label{sec:diamond}

We review certain structure results related to our dual pairs \cite{Howe.1989.transcending}.

Let $ W_{\R} \simeq \R^{2 N} $ be a real symplectic space which realizes $ \bbG = Sp(2N,\R ) $ as a symplectic group on $ W_{\R} $.  
There is a canonical complex structure on $ W_{\R} $ and we can view $ W_{\R} $ as (the underlying real vector space of) a complex vector space $ W \simeq \C^N $.  
The symplectic form on $ W_{\R} $ is then given by the imaginary part of a canonical positive definite Hermitian form on $ W $.  
By this identification, a maximal compact subgroup $ \bbK $ of $ \bbG $ is realized as the unitary group $ U( W ) \simeq U( N ) $ on $ W \simeq \C^N $.

We may choose maximal compact subgroups $ \Kbig $ and $ \Ksmall $ of $ \Gbig $ and $ \Gsmall $ respectively, in such a way that $ \Kbig \cdot \Ksmall $ is contained in the standard maximal compact subgroup $ \bbK \simeq U(N) $ of $ \bbG $.  

In view of Table \ref{table:dual.pairs.treated}, we sometimes write $ G = G(p, q) $ for $ O(p, q), U(p, q) \text{ or } Sp(p, q) $.  
We denote simply $ G(k) = G(k, 0) $, which is compact. 

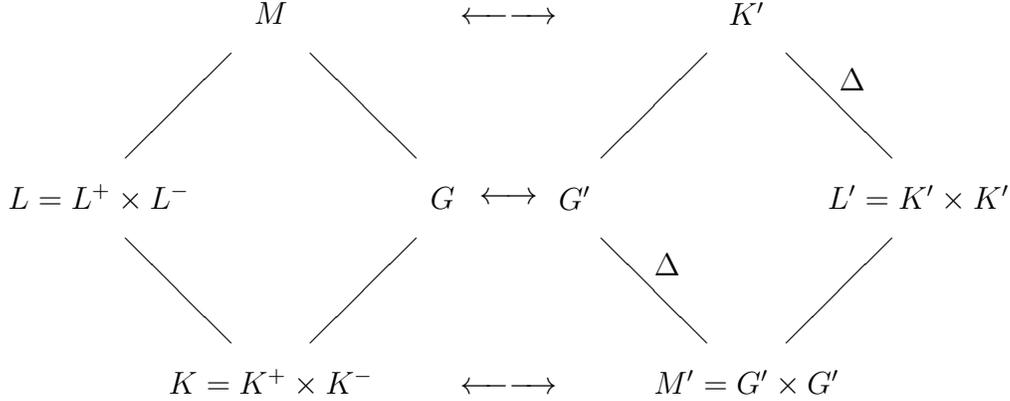
\begin{figure}[htbp]
\caption{The diagram of a diamond pair.}
\label{fig:diamond.pair}
\hfil
\begin{picture}(400,190)

\put( 55, 115){\line(1, 1){40}}
\put( 55,  85){\line(1,-1){40}}
\put(125, 155){\line(1,-1){40}}
\put(125,  45){\line(1, 1){40}}

\put(235, 115){\line(1, 1){40}}
\put(235,  85){\line(1,-1){40}}
\put(305, 155){\line(1,-1){40}}
\put(305,  45){\line(1, 1){40}}

\put(330, 145){\makebox(0,0){$ \Delta $}}
\put(260,  75){\makebox(0,0){$ \Delta $}}

\put( 45, 100){\makebox(0,0){$ \Lbig = \Lbig^{+} \times \Lbig^{-} $}}
\put(110, 170){\makebox(0,0){$ \Mbig $}}
\put(110,  30){\makebox(0,0){$ \Kbig = \Kbig^{+} \times \Kbig^{-} $}}
\put(175, 100){\makebox(0,0){$ \Gbig $}}

\put(225, 100){\makebox(0,0){$ \Gsmall $}}
\put(290, 170){\makebox(0,0){$ \Ksmall $}}
\put(290,  30){\makebox(0,0){$ \Msmall = \Gsmall \times \Gsmall $}}
\put(355, 100){\makebox(0,0){$ \Lsmall = \Ksmall \times \Ksmall $}}

\put(200,  30){\makebox(0,0){$ {\xleftarrow{\quad}}{\xrightarrow{\quad}} $}}
\put(200, 100){\makebox(0,0){$ \longleftrightarrow $}}
\put(200, 170){\makebox(0,0){$ {\xleftarrow{\quad}}{\xrightarrow{\quad}} $}}

\end{picture}
\hfil
\end{figure}

From the pair $ ( \Gbig, \Gsmall ) $, we define another three dual pairs, 
which form so-called \textit{diamond dual pairs} (see \cite[\S 5]{Howe.1989.transcending}).  
Namely, take the commutant of $ \Kbig $ in $ \bbG $ and denote it by $ \Msmall $.  
Then $ \Msmall $ is also of Hermitian symmetric type and isomorphic to $ \Gsmall \times \Gsmall $ containing $ \Gsmall $ as the diagonal.  
The pair $ ( \Kbig , \Msmall ) $ is a dual pair of compact type.  
Also, take the full commutant of $ \Kbig $ in $ \bbK $, and denote it by $ \Lsmall $.  
Then $ \Lsmall $ is a maximal compact subgroup of $ \Msmall $,  
and it is isomorphic to $ \Ksmall \times \Ksmall $ and contains $ \Ksmall $ as a diagonal subgroup.  
Similarly, define $ \Mbig $ as the commutant of $ \Ksmall $ in $ \bbG $, 
and $ \Lbig $ the commutant in $ \bbK $.  
Let us summarize this somewhat complicated situation by the following diagram (Fig.~\ref{fig:diamond.pair}).  
An explicit description of the diamond pairs for our three cases is given in  \cite[(5.3)]{Howe.1989.transcending}.  
In the table there, $ \Lbig $ (resp.\  $ \Lsmall $) in our notation is written as $ M^{(1, 1)} $ (resp.\  $ M'{}^{(1,1)} $).  
For convenience of readers, we reproduce the table here (Table \ref{table:diamond.dual.pairs} below) with some additional remarks.
\begin{table}[hbtp]
\caption{Three diamond dual pairs.}
\label{table:diamond.dual.pairs}
\begin{gather*}
\text{\footnotesize $
\begin{array}{@{}ccccc@{}}
( \Gbig , \Gsmall ) & \Mbig & \Kbig = \Kbig^{+} \times \Kbig^{-} & \Lbig = \Lbig^{+} \times \Lbig^{-} & \Ksmall \\
\hline
\\
( O( p, q ) , Sp( 2 n, \R ) ) & U(p, q) & O(p) {\times} O(q) & U(p) {\times} U(q) & U(n) \\[5pt]
( U( p, q ), U( m, n ) ) ) & U(p, q) {\times} U(p, q) & U(p) {\times} U(q) & ( U(p) {\times} U(p) ) {\times} ( U(q) {\times} U(q) ) & U(m) {\times} U(n) \\[5pt]
( Sp( p, q ) , O^{\ast}( 2 n ) ) & U( 2 p , 2 q ) & Sp(p) {\times} Sp(q) & U( 2 p ) {\times} U( 2 q ) & U(n) 
\end{array}
$}\\[10pt]
\Msmall = \Gsmall \times \Gsmall , \hspace*{.1\textwidth}  \Lsmall = \Ksmall \times \Ksmall 
\end{gather*}
\end{table}

Note that $ \Mbig / \Lbig $ is a Hermitian symmetric space (not necessarily irreducible), 
and the containment $ \Mbig \supset \Gbig $ (resp.{} $ \Lbig \supset \Kbig $) may or may not be the diagonal map, but it is always a symmetric pair.

Recall $ \bbG = Sp( W_{\R} ) $ and the complex vector space $ W \simeq \C^N $ which is identical to $ W_{\R} $ as a real vector space.  
Then there exists a direct sum decomposition $ W_{\R} = W^+_{\R}\oplus W^-_{\R} $
and correspondingly $ W = W^+ \oplus W^- $, which are compatible
with direct product decompositions 
$ \Lbig = \Lbig^+ \times \Lbig^- $ and 
$ \Kbig = \Kbig^+ \times \Kbig^- $ in the following way.  
The subgroups $ \Lbig^{\pm} $ and $ \Kbig^{\pm} $ are contained in the unitary group $ U( W^{\pm} ) $.  
Moreover the pairs $ ( \Kbig^{\pm} , \Gsmall ) $ and $ ( \Lbig^{\pm} , \Ksmall ) $ are dual pairs in $ Sp( W^{\pm}_{\R} ) $.  
Also note that $ \Lbig^{\pm} \supset \Kbig^{\pm} $ is a symmetric pair.  

For explicit realization of these decompositions, see the Appendix.

\section{Theta lift of a nilpotent orbit}
\label{sec:lifting.nilpotent.orbit}

Let $ \liec{G} = \Lie ( \bbG )_{\C} $ (resp.\  $ \liec{K} = \Lie ( \bbK )_{\C} $) be 
the complexified Lie algebra of $ \bbG = Sp(2N, \R) $ (resp.\  $ \bbK = U(N) $).  
Let $ \liec{G} = \liec{K} \oplus \liec{P} $ be the corresponding (complexified) Cartan decomposition.  
Note that the complexified Lie group $ \bbKC $ acts on the Cartan space $ \liec{P} $ by the restriction of the adjoint action; 
under this action $\liec{P} $ breaks up into $\lie{P}_{\pm}$, 
each of which can be identified with a copy of the space $ \Sym_N(\C) $ of complex symmetric matrices of size $ N $. 
The action of $ k \in GL(N, \C) \simeq \bbKC $ on 
$ \liec{P} = \lie{P}_{+} \oplus \lie{P}_{-} \simeq \Sym_N(\C) \oplus \Sym_N(\C) $ is given by 
\begin{displaymath}
k \cdot ( X , Y ) = (  k X \transpose{k} , \transpose{k}^{-1}Y k^{-1}  ) \qquad ( X , Y \in \Sym_N(\C) ) .
\end{displaymath}
%

%
%

Let $ \Omin \subseteq \lie{P}_{+}$ be the minimal nilpotent $ \bbKC \simeq GL(N, \C) $ orbit 
in $ \lie{P}_+ $. We have
the identification 
\begin{math}
\Omin = \{  (X,0) \mid X \in \Sym_N(\C) , \linebreak[3] \rank X = 1 \} .
\end{math}
It is also well-known that there is a map 
\begin{displaymath}
\begin{array}{rccc}
j : \; & W & \longrightarrow  & \closure{\Omin} \\
       & \vertin & & \vertin \\
       & w &  \longmapsto & (\transpose{w} w, 0 ) 
\end{array} 
\qquad w = ( w_1, \ldots, w_N ) \in \C^N = W
\end{displaymath}
by which $ \closure{\Omin} = \Omin \cup \{ 0 \} $ is regarded as a geometric quotient of the space $ W $ by the action of 
$ O(1,\C) = \{ \pm 1 \} $,  
and that the quotient map is $ \bbK_{\C} $-equivariant.  

Let us consider the dual pair $ ( \Gbig , \Gsmall ) \subseteq \bbG $.  
The maximal compact subgroup $ \Kbig $ (resp.~$ \Ksmall $) determines 
a Cartan decomposition of the complexified Lie algebra $ \gbig $ (resp.~$ \gsmall $) of $ \Gbig $ (resp.~$ \Gsmall $): 
\begin{displaymath}
\gbig   = \kbig   \oplus \sbig   , \qquad 
\gsmall = \ksmall \oplus \ssmall , 
\end{displaymath}
where $ \kbig $ (resp.~$ \ksmall $) is the complexified Lie algebra of $ \Kbig  $ 
(resp.~$ \Ksmall  $), 
and $ \sbig $ (resp.~$ \ssmall $) is the orthogonal complement of $ \kbig $ (resp.~$ \ksmall $) 
with respect to the Killing form.  
The subspace $ \sbig $ (resp.~$ \ssmall $) carries a linear action of $ \KbigC $ (resp.~$ \KsmallC $) 
via the restriction of the adjoint representation.  
Since $ \Gsmall / \Ksmall $ is a Hermitian symmetric space by assumption, 
$ \ssmall $ breaks up into irreducible pieces $ \ssmall_{\pm} $ under the action of $ \KsmallC $.   

We can arrange the Cartan decompositions so that $ \sbig \oplus \ssmall $ is contained in $ \liec{P} $.  
Then the inclusion map $ \sbig \subset \liec{P} $ (resp.{} $ \ssmall \subset \liec{P} $) induces 
the projection $ \liec{P}^{\ast} \rightarrow \sbig^{\ast} $ (resp.{} $ \liec{P}^{\ast} \rightarrow \ssmall^{\ast} $).  
By identifying  the complex linear dual $ \liec{P}^{\ast} $ (resp. $ \sbig^{\ast}$, 
$ \ssmall^{\ast}$)
with $ \liec{P} $ (resp. $ \sbig $, $ \ssmall$) by the invariant trace form of matrices, 
we obtain two projection maps $P_{\sbig}: \liec{P} \rightarrow \sbig $ and
$P_{\ssmall}: \liec{P} \rightarrow \ssmall $.
The projections restricted to $ \closure{\Omin} $ will then induce the so-called moment maps (cf.~\cite{Ohta.preprint1}):  
\begin{align*}
\momentbig : \ \ & W \xrightarrow{\begin{subarray}{c}
\text{geom. quotient}\\
\text{by $ \{ \pm 1 \} $}
\end{subarray}
} \closure{\Omin} \xrightarrow{\quad\text{projection}\quad} \sbig \\[3pt]
\momentsmall : \ \ & W \xrightarrow{\hspace{.125\textwidth}} \closure{\Omin} \xrightarrow{\hspace{.13\textwidth}} \ssmall 
\end{align*}
These maps are $ \KbigC \times \KsmallC $-equivariant with the trivial $ \KbigC $-action on $ \ssmall $, and 
the trivial $ \KsmallC $-action on $ \sbig $ respectively. Also we can (and we will) choose the subspaces $ \ssmall_{\pm} \subset \ssmall $ and $ W^{\pm} \subset W $ 
so that 
\begin{math}
\momentsmall ( W^{\pm} ) \subset \ssmall_{\pm}  
\end{math}
holds (cf.~\S \ref{sec:diamond} for the notation $ W^{\pm} $).  

For a set $ S \subseteq \liec{P} $, we denote by $ \nilpotents( S ) $ the subset of all nilpotent elements in $ S $.  

\begin{lemma}
We can choose Cartan decompositions and moment maps so that the following conditions hold.  
\begin{thmenumerate}
\item
$ \momentsmall $ is an affine quotient map by $ \KbigC $ onto $ \ssmall $, i.e., 
$ \ssmall \simeq W \GITquotient \KbigC $.  
\item
$ \momentbig $ is an affine quotient map by $ \KsmallC $ onto its closed image.
\item
$ \momentbig( \momentsmall^{-1}( \nilpotents( \ssmall ) ) ) \subset \nilpotents( \sbig ) $.
\end{thmenumerate}
\end{lemma}

\begin{remark}
For part (1) to hold true, we only need to have the following conditions:  $\min (p, q) \geq n $ for $( O( p, q ) , Sp( 2 n, \R ) )$,   $\min (p, q)\geq \min (m, n)$ for $( U( p, q ), U( m, n ) )$, and $\min (p,q) \geq \frac {n}{2}$ for $(Sp( p, q ) , O^{\ast}( 2 n ) )$. They are much weaker than the stable range condition. 
\end{remark}

\begin{proof}
To show that $ \momentsmall : W \rightarrow \ssmall $ is an affine quotient map, 
it suffices to show that $ \momentsmall $ induces an algebra isomorphism 
$ \C[ \ssmall ] \simeq \C[ W ]^{\KbigC} $, where $ \C[ W ]^{\KbigC} $ denotes the space
of invariants of 
the polynomial ring $ \C[ W ] $. Clearly it suffices to show that the restrictions $ \momentsmall \restrict_{W^{\pm}} $ induce algebra isomorphisms 
$ \C[ \ssmall _{\pm}] \simeq \C[ W^{\pm} ]^{\KbigpqC} $.

Similarly, statement (2) asserts that the induced algebra homomorphism 
$ \momentbig^{\ast} : \C[ \sbig ] \rightarrow \C[ W ] $ maps 
$ \C[ \sbig ] $ onto the $ \KsmallC $-invariants $ \C[ W ]^{\KsmallC} $.  
Note that the closed subvariety $ \Im \momentbig \subset \sbig $ is then defined by the ideal of 
relations of generators of invariants, 
or in other words, the kernel of $ \momentbig^{\ast} $.  

Explicit constructions of moment maps for each case of the three dual pairs in Table \ref{table:dual.pairs.treated} 
are summarized in the Appendix.  
From the explicit formulas given there, the statements (1) -- (3) follow easily 
from classical invariant theory \cite{Weyl.1946}.
\end{proof}

Using the moment maps, we shall define the theta lift of nilpotent orbits 
from $ \nilpotents( \ssmall ) $ to $ \nilpotents( \sbig ) $.
Before that, we review the following well-known result due to Kostant.

We consider the null cones
\begin{equation}
\nullcone = \momentsmall^{-1}( 0 ) \subset W , \qquad 
\nullconepq = ( \momentsmallpq )^{-1}( 0 ) \subset \Wpq , 
\end{equation}
where $ \momentsmallpq : \Wpq \rightarrow \ssmallpq $ is the restriction of $ \momentsmall $ to $ \Wpq $.
Let $ J^{\pm} \subset \C[ W^{\pm} ] $ be the augmentation ideal 
generated by homogeneous invariants $ \C[ W^{\pm} ]_{+}^{ \KbigpqC } $ of positive degree, 
and let $ J = \C[ W ]_{+}^{ \KbigC } \cdot \C[W] = J^+\otimes J^{-}$.  

We observe that the stable range condition implies that
\begin{displaymath}
\text{
$ {J}^{\pm}$ is the ideal in $ \C[W^{\pm}] $ vanishing on $ \nullconepq $ and $ {J}^{\pm} $ is prime;}
\end{displaymath}
\begin{displaymath}
\text{
The single $\KbigpqC$-orbit of generic elements of $ \nullconepq $ is dense in $ \nullconepq $.}
\end{displaymath}
See for example Proposition 3.7.3.7, and Theorem 3.8.6.4 in \cite{Howe.1995.SchurLecture}, and references therein.  
Some relevant properties of the null cones may also be found in the Appendix. 

As usual, we denote by $ \harmonics( \KbigpqC ) $ the space of harmonic polynomials on $ \Wpq $ 
{\upshape (}under the action of $\KbigpqC ${\upshape )}. Similar notations apply throughout
this article. Applying Kostant's criterion  \cite{Kostant.1963}, we get

\begin{theorem}
\label{thm:harmonics.and.invariants}
Let $ ( G, G' ) $ be in the stable range $($see {\upshape Table \ref{table:dual.pairs.treated}}$)$ with $ G' $ the smaller member, 
and consider the dual pairs $ ( \Kbigpq , \Gsmall ) $ in $ Sp( \Wpq_{\R} ) $.  
Then with respect to the action of $ \KbigpqC $, we have a tensor product decomposition of $ \C[ \Wpq ] $:
\begin{displaymath}
\C[ \Wpq ] \simeq \harmonics( \KbigpqC ) \otimes \C[ \Wpq ]^{\KbigpqC}.
\end{displaymath}
\end{theorem}

\begin{corollary}
\label{cor:moment.map.is.flat}
We have
\begin{displaymath}
\C[ W ] \simeq \harmonics( \KbigC ) \otimes \C[ W ]^{ \KbigC }
\simeq \harmonics( \KbigC ) \otimes \C[ \ssmall ] .
\end{displaymath}
Consequently the moment map $ \momentsmall : W \rightarrow \ssmall $ is flat.
\end{corollary}

\begin{proof}
Note that the induced map $\momentsmall^{\ast}$ of the moment map $ \momentsmall $ is just an inclusion: 
\begin{math}
\momentsmall^{\ast} : \C[ \ssmall ] \simeq \C[ W ]^{ \KbigC } \hookrightarrow \C[ W ] .
\end{math}
Since $\C[ W ] \simeq \harmonics( \KbigC ) \otimes \C[ \ssmall ]$, it is free over $ \C[ \ssmall ] $.   
This implies that $ \momentsmall $ is flat.
\end{proof}

We are now in a position to prove

\begin{theorem}
\label{thm:fiber.of.Osmall}
Take a nilpotent $ \KsmallC $-orbit $ \Osmall \subset \nilpotents( \ssmall ) $.  
Then the scheme theoretic fiber 
$ \momentsmall^{-1}( \closure{\Osmall} ) = W \times_{ \ssmall } \closure{ \Osmall } $ 
is a reduced, closed irreducible affine subvariety of $ W $.  
Moreover, it is the closure of a single $ \KbigC \times \KsmallC $-orbit in $ W $.  
\end{theorem}

\begin{proof}
Let us prove that the fiber 
$ W \times_{ \ssmall } \closure{ \Osmall } $ is reduced.  
For this, it is enough to show that the fiber of each closed point in $ \ssmall $ is reduced.  
Namely, we prove the following lemma.

\begin{lemma}
Let $ x \in \ssmall $ be a closed point.  
Then the scheme theoretic fiber 
\begin{math}
\momentsmall^{-1}( x ) = W \times_{ \ssmall } \{ x \} 
\end{math}
is a closed, reduced and irreducible affine subvariety of $ W $.  
Moreover, it is the closure of a single $ \KbigC $-orbit.  
\end{lemma}

\begin{proof}
The scheme theoretic fiber is defined by 
\begin{math}
W \times_{\ssmall} \{ x \} = \Spec \left( \C[ W ] \otimes_{ \C[ \ssmall ] } \C_x \right) .  
\end{math}
Here $ \C_x = \C[ \ssmall ] / \mathfrak{m}_x \simeq \C $, 
where $ \lie{m}_x $ denotes the maximal ideal corresponding to the closed point $ x $.  
Put $ \mathcal{A}_x = \C[ W ] \otimes_{ \C[ \ssmall ] } \C_x $.  
We are to show that $ \mathcal{A}_x $ is an integral domain.  
The proof is based on the general argument about deformations of homogeneous integral domains.  

Let us recall the tensor product decomposition 
\begin{math}
\C[ W ] = \harmonics( \KbigC ) \otimes \C[ \ssmall ].
\end{math}
Thus, if we abbreviate $ \harmonics = \harmonics( \KbigC ) $, we can identify 
\begin{displaymath}
\mathcal{A}_x = \C[ W ] \otimes_{ \C[ \ssmall ] } \C_x = \harmonics \otimes \C_x \simeq \harmonics .
\end{displaymath}
However, note that $ \harmonics $ is not an algebra 
but it only enjoys a structure of $ \KbigC \times \KsmallC $-module.  
Since $ \harmonics $ is a graded subspace of $ \C[ W ] $, 
it is naturally graded by the ordinary degree of polynomials.  
We denote the grading by 
\begin{math}
\harmonics \simeq \bigoplus\nolimits_{ k \geq 0 } \harmonics_k .
\end{math}
Then, $ \mathcal{A}_x $ becomes a filtered algebra by putting 
\begin{displaymath}
( \mathcal{A}_x )_i = \sum\nolimits_{ k \leq i } \harmonics_k \otimes 1 \subset \harmonics \otimes \C_x = \mathcal{A}_x .
\end{displaymath}
\indent
Now assume that $ \mathcal{A}_x $ is not an integral domain.  
Then there are non-zero elements $ a, b \in \mathcal{A}_x $ such that $ a \cdot b = 0 $.  
Let us express $ a $ and $ b $ as 
\begin{displaymath}
\textstyle 
a = \sum\nolimits_{ i = 0 }^d a_i \otimes 1 \in \harmonics \otimes \C_x \; ; \quad 
b = \sum\nolimits_{ j = 0 }^{d'} b_j \otimes 1 \in \harmonics \otimes \C_x ,
\end{displaymath}
with $ a_d \neq 0 $ and $ b_{d'} \neq 0 $.  
Then $ a b = 0 $ implies 
\begin{math}
( a_d b_{d'} ) \otimes 1 \in ( \mathcal{A}_x )_{ d + d' - 1 } ,
\end{math}
but this is impossible.  
To see it, observe that the null cone $ \nullcone = \nullconep \times \nullconeq $ is defined by the augmentation ideal $J$ 
which is prime, and so $ \nullcone $ is irreducible.  
Note that $ \C[ \nullcone ] = \C[ W ] / J \simeq \harmonics $.  
Since $ a_d b_{d'} \not\in \harmonics $ by the irreducibility of the null cone, 
we have 
\begin{displaymath}
\textstyle
a_d \cdot b_{d'} \otimes 1 \not\in \sum_{k = 0}^{d + d' - 1} \harmonics_k \otimes \C_x = ( \mathcal{A}_x )_{d + d' - 1} ,
\end{displaymath}
which is a contradiction.  

\newcommand{\orbittemp}{\mathcal{O}}

Next, we prove that the fiber $ \momentsmall^{-1}(x) $ contains an open dense $ \KbigC $-orbit.  
Put $ M = \momentsmall^{-1}(x) $ and denote by $ \widehat{M} $ the asymptotic cone of $ M $ 
(see \cite[\S 5.2]{Popov.Vinberg.1994} for the definition of asymptotic cone).  
Then, by the flatness of $ \momentsmall $, the asymptotic cone $ \widehat{M} $ coincides with the null cone $ \nullcone $.  
Let $ \orbittemp $ be a generic $ \KbigC $-orbit in $ M $.  
Consider the cone $ \C M $ generated by $ M $ in 
$ W $, then it is clear that the dimension of a generic orbit in $ \C M $ is equal to $ \dim \orbittemp $, which  
in turn coincides with the dimension of the generic orbit in $ \closure{\C M} \subset W $.  
Since $ \nullcone = \widehat{M} \subset \closure{\C M} $, the dimension of a generic orbit in $ \nullcone $ cannot exceed that of $ \orbittemp $.  
Note that $ \nullcone $ has an open dense orbit (see the Appendix for details).  
This means that $ \dim \orbittemp \geq \dim \nullcone $.  
On the other hand, we have an equality $ \dim M = \dim \widehat{M} = \dim \nullcone $ of dimensions, 
hence $ \dim \orbittemp \geq \dim M $.  
Since $ \orbittemp \subset M $, we conclude that $ \dim \orbittemp = \dim M $, and 
that $ \orbittemp $ is an open dense orbit in $ M $, by the irreducibility of $ M $ just proved above.  
\end{proof}

Let us return to the proof of Theorem \ref{thm:fiber.of.Osmall}.  

By $ \KsmallC $-equivariancy of $ \momentsmall $, we get 
\begin{math}
\momentsmall^{-1}( \Osmall ) = \KsmallC \cdot \momentsmall^{-1}( \{ x \} ) 
\end{math}
for any $ x \in \Osmall $.
Consider the multiplication map 
\begin{math}
\KsmallC \times \momentsmall^{-1}( \{ x \} ) \rightarrow \momentsmall^{-1}( \Osmall ) .
\end{math}
Since $ \KsmallC \times \momentsmall^{-1}( \{ x \} ) $ is irreducible by the above lemma, 
$ \momentsmall^{-1}( \Osmall ) $ is also irreducible (as an image of an irreducible set).  
Therefore its closure 
$ \closure{ \momentsmall^{-1}( \Osmall ) } $ 
is irreducible.
Since the moment map $ \momentsmall $ is flat by \corollaryref{cor:moment.map.is.flat}, 
it is an open map (\cite[Ex.~(III.9.1)]{Hartshorne.1977}).  
Thus we conclude that $ \momentsmall^{-1}( \closure{ \Osmall } ) = \closure{ \momentsmall^{-1}( \Osmall ) } $ 
is irreducible.  

\newcommand{\orbittemp}{\mathcal{O}}

By the same lemma, 
$ \momentsmall^{-1}( x ) $ contains an open dense $ \KbigC $-orbit $ \orbittemp_x $.  
Then the union of $ \KsmallC $-translates $ \KsmallC \cdot \orbittemp_x $ is dense in $ \momentsmall^{-1}( \Osmall ) $.  
Thus 
$ \closure{ \KsmallC \cdot \orbittemp_x } = \closure{ \momentsmall^{-1}( \Osmall ) } 
	= \momentsmall^{-1}( \closure{ \Osmall } ) $, which means that 
$ \KsmallC \cdot \orbittemp_x $ is open in $ \momentsmall^{-1}( \closure{ \Osmall } ) $.
\end{proof}

\begin{remark}
As mentioned in the Introduction, T.~Ohta proved the irreducibility and the existence of an open dense orbit by a totally different (but case-by-case) method 
(\cite{Ohta.preprint1}, \cite{Ohta.preprint2}).  
In fact, his method is applicable beyond the stable range.  
However, outside the stable range, 
things get much more complicated and the same statement of the above theorem is no longer true.  
\end{remark}

From the above theorem, we see that 
$ \momentbig( \momentsmall^{-1} ( \closure{ \Osmall } ) ) $ is a $ \KbigC $-stable irreducible closed set in $ \sbig $. Since $ \momentbig( \momentsmall^{-1} ( \closure{ \Osmall } ) )$ is contained in the nilpotent variety $ \nilpotents( \sbig ) $, which has only finite number of $ \KbigC $-orbits, it must be the closure of a single $ \KbigC $-orbit $ \Obig $.

\begin{definition}
Let $ \Osmall \subset \ssmall $ be a nilpotent $ \KsmallC $-orbit.  
Then $ \momentsmall^{-1}( \closure{\Osmall} ) $ is a closed, irreducible affine subvariety in $ W $.  
The image $ \momentbig( \momentsmall^{-1}( \closure{\Osmall} ) ) $ is 
the closure of a single nilpotent $ \KbigC $-orbit in $ \sbig $, i.e., 
$ \momentbig( \momentsmall^{-1}( \closure{\Osmall} ) ) = \closure{ \Obig } $ for a certain nilpotent 
$ \KbigC $-orbit $ \Obig $.  
We call $ \Obig $ the \textit{theta lift} of $ \Osmall $, and denote it by $ \Thetalift( \Osmall ) $.
\end{definition}

We introduce some notations. For any complex reductive Lie group $E$, let $\irreps{E}$ be the set of equivalence classes of irreducible finite dimensional representations of $E$. Given
a (completely reducible) representation of $E$ on $\mathcal U$, let $\irreps{ E; \mathcal U}$  be the subset of $\irreps{E}$ which appear in $ \mathcal U$. 


%
%
%
%

We now recall some standard results of Howe \cite{Howe.1989.classical}. As noted in the Introduction, the space of harmonics $ \harmonics( \KbigpqC ) $ is multiplicity-free as a representation of
$\KbigpqC \times \KsmallC$. Furthermore 
the decomposition 
\begin{equation}
\label{eqn:decomp.harmonics.by.KtimesKprime}
\harmonics( \KbigpC )\restrict _{\KbigpC \times \KsmallC} \simeq 
	\directsum\nolimits_{ \sigma \in \irreps{ \KbigpC; \harmonics^{+} }, 
		\tau \in \irreps{ \KsmallC ; \harmonics^{+} } 
		} \sigma \boxtimes \tau
\end{equation}
determines a one-to-one correspondence between 
$ \sigma \in \irreps{ \KbigpC ; \harmonics^{+} } $ 
and $ \tau \in \irreps{ \KsmallC ; \harmonics^{+} } $, 
where we abbreviate $ \harmonics( \KbigpC ) $ to $ \harmonics^{+} $.  
Similar statements apply to $ \harmonics( \KbigqC ) $. 

We shall abbreviate $ R^{\pm}( \KsmallC ) = \irreps{ \KsmallC ; 
\harmonics( \KbigpqC ) } $ 
in the following. Under the assumption of stable range, 
one can explicitly check that $ \tau^{\ast} \in R^{+}( \KsmallC ) $ if and only if 
$ \tau \in R^{-}( \KsmallC ) $, 
where  $ \tau^{\ast} $ denotes the contragredient representation of $ \tau \in \irreps{K'} $.  
So we put 
\begin{equation}
\label{eqn:definition.R}
R( \KsmallC ) = R^{+}( \KsmallC )^{\ast} = R^{-}( \KsmallC ) .
\end{equation}
To summarize, for each $ \tau \in R^{\pm}( \KsmallC ) $, 
there exists a unique $ \sigma \in \irreps{\KbigpqC} $ such that 
\begin{displaymath}
\Hom_{\KbigpqC \times \KsmallC}( \sigma \boxtimes \tau, \harmonics( \KbigpqC ) ) \neq 0 .
\end{displaymath}
We denote this $ \sigma $ by $ \sigma^{\pm}( \tau ) $, specifying the dependency of $ \tau $ 
and the sign $ \pm $.  
Then, we can rewrite \eqnref{eqn:decomp.harmonics.by.KtimesKprime} and its sister statement as 
\begin{equation}
\label{eqn:rewrited.decomp.harmonics.by.KtimesKprime}
\harmonics( \KbigpC ) 
\simeq \directsum_{ \tau \in R( \KsmallC ) } \sigma^{+}( \tau^{\ast} ) \boxtimes \tau^{\ast} , 
\quad \text{ and } \quad 
\harmonics( \KbigqC ) 
\simeq \directsum_{ \tau \in R( \KsmallC ) } \sigma^{-}( \tau ) \boxtimes \tau .
\end{equation}

\begin{remark}
By the theory of highest weight, we may identify $ R( \KsmallC ) $ with a
certain semigroup in the integral weight lattice of $ \KsmallC $.  We note
that in an explicit realization, the semigroup $ R( \KsmallC ) $ may be
identified with a set of partitions, and it is saturated. The explicit
decompositions in Equation
(\ref{eqn:rewrited.decomp.harmonics.by.KtimesKprime}) are well-known. See
\cite{Kashiwara.Vergne.1978} for the orthogonal and unitary cases, or
\cite{Howe.1995.SchurLecture} in general.
\end{remark}

\begin{theorem}
\label{thm:fun.ring.of.lifted.orbit}
Let $ ( G, G' ) $ be in the stable range 
$($see {\upshape Table \ref{table:dual.pairs.treated}}$)$ with $ G' $ the smaller member.  
For a nilpotent $ \KsmallC $-orbit $ \Osmall \subset \nilpotents( \ssmall ) $, denote its theta lift by 
$ \Obig = \theta( \Osmall ) $.  Put $ \Xi( \Osmall ) = \momentsmall^{-1}( \closure{\Osmall} ) $. 
Then the closure $ \closure{ \Obig } $ is an affine quotient of $ \Xi( \Osmall ) $ by $ \KsmallC $.  
\begin{displaymath}
\closure{\Obig} \simeq \Xi( \Osmall ) \GITquotient \KsmallC .
\end{displaymath}
Moreover, the $ \KbigC $-module structure of the regular function ring of $ \closure{ \Obig } $ 
is given in terms of $ \closure{ \Osmall } $ as follows:  
\begin{align}
\C[ \closure{ \Obig } ] 
&\simeq \left( \harmonics( \KbigC ) \otimes \C[ \closure{ \Osmall } ] \right)^{\KsmallC} \nonumber \\
&\simeq \directsum_{ \tau_1, \tau_2 \in R( \KsmallC ) } 
\Hom_{\KsmallC} ( \tau_1 \otimes \tau_2^{\ast} , \C[ \closure{ \Osmall } ] ) 
	\otimes \bigl( \sigma^{+}( \tau_1^{\ast} ) \boxtimes \sigma^{-}( \tau_2  ) \bigr) , 
\end{align}
where $ \sigma^{+}( \tau_1^{\ast} ) \boxtimes \sigma^{-}( \tau_2 ) $ is an irreducible representation of 
$ \KbigC = \KbigpC \times \KbigqC $ given in Equation {\upshape \eqnref{eqn:rewrited.decomp.harmonics.by.KtimesKprime}}, 
and $ \KbigC $ acts on the space of multiplicities 
$ \Hom_{\KsmallC} ( \tau_1 \otimes \tau_2^{\ast} , \C[ \closure{ \Osmall } ] )$ 
trivially.
\end{theorem}

\begin{proof}
 The assertion that $ \closure{ \Obig } $ is an affine quotient follows from general theory of affine quotient maps 
because $ \momentbig $ itself is a quotient map, and 
$\Xi( \Osmall ) $ is a $ \KsmallC $-invariant, affine closed subvariety of $ W $  
(cf.~\cite[Prop.~3.3]{Nishiyama.JSPS.DFG.proc}).

Next we prove the statement on the module structure of $ \C[ \closure{ \Obig } ] $.  
By \theoremref{thm:fiber.of.Osmall}, we have 
\begin{align}
\C[ \Xi( \Osmall ) ] & \simeq \C[ W \times_{ \ssmall } \closure{ \Osmall } ] 
	= \C[ W ] \otimes_{ \C[ \ssmall ] } \C[ \closure{ \Osmall } ] \notag \\
	& \simeq \left( \C[ \Wp ] \otimes \C[ \Wq ] \right) \otimes_{ \C[ \ssmallp ] \otimes \C[ \ssmallq ] } \C[ \closure{ \Osmall } ] \notag \\
	& \simeq \left( \C[ \Wp ] \otimes_{ \C[ \ssmallp ] } \C[ \closure{ \Osmall } ] \right) 
		\otimes_{ \C[ \closure{ \Osmall } ] } 
		\left( \C[ \Wq ] \otimes_{ \C[ \ssmallq ] } \C[ \closure{ \Osmall } ] \right) .
\label{eqn:Wp.tensor.Wq}
\end{align}
Note that 
\begin{math}
\C[ \Wpq ] \simeq \harmonics( \KbigpqC ) \otimes \C[ \Wpq ]^{\KbigpqC} 
\end{math}
and $ \C[ \Wpq ]^{\KbigpqC} \simeq \C[ \ssmallpq ] $ by 
\theoremref{thm:harmonics.and.invariants}.  
Therefore we have 
\begin{displaymath}
\C[ \Wp ] \otimes_{ \C[ \ssmallp ] } \C[ \closure{ \Osmall } ] 
	\simeq \bigl( \harmonics( \KbigpC ) \otimes \C[ \Wp ]^{\KbigpC} \bigr) 
		\otimes_{ \C[ \ssmallp ] } \C[ \closure{ \Osmall } ] 
	\simeq \harmonics( \KbigpC ) \otimes \C[ \closure{ \Osmall } ] , 
\end{displaymath}
and \eqnref{eqn:Wp.tensor.Wq} becomes
\begin{displaymath}
\left( \harmonics( \KbigpC ) \otimes \C[ \closure{ \Osmall } ] \right) \otimes_{ \C[ \closure{ \Osmall } ] } 
\left( \harmonics( \KbigqC ) \otimes \C[ \closure{ \Osmall } ] \right) 
\simeq \left( \harmonics( \KbigpC ) \otimes \harmonics( \KbigqC ) \right) \otimes \C[ \closure{ \Osmall } ] .
\end{displaymath}
Since $ \momentbig $ is an affine quotient map by the action of $ \KsmallC $, we have 
\begin{align*}
\C[ \closure{ \Obig } ] & \simeq \C[ \Xi( \Osmall ) ]^{\KsmallC} 
		\simeq  \left( \left( \harmonics( \KbigpC ) \otimes \harmonics( \KbigqC ) \right) 
		\otimes \C[ \closure{ \Osmall } ] \right)^{\KsmallC} \\
	&\simeq \directsum_{ \tau_1, \tau_2 \in R( \KsmallC ) } 
		\left( \bigl( \sigma^{+}( \tau_1^{\ast} ) \boxtimes \tau_1^{\ast} \bigr) 
		\otimes \bigl( \sigma^{-}( \tau_2 ) \boxtimes \tau_2 \bigr) 
		\otimes \C[ \closure{ \Osmall } ] \right)^{\KsmallC} \\
	&\simeq \directsum_{ \tau_1, \tau_2 \in R( \KsmallC ) } 
		\left( \tau_1^{\ast} \otimes \tau_2  \otimes \C[ \closure{ \Osmall } ] \right)^{\KsmallC} 
		\otimes \bigl( \sigma^{+}( \tau_1^{\ast} ) \boxtimes \sigma^{-}( \tau_2 ) \bigr) \\
	&\simeq \directsum_{ \tau_1 , \tau_2 \in R( \KsmallC ) } 
		\Hom_{\KsmallC} ( \tau_1 \otimes \tau_2^{\ast} , \C[ \closure{ \Osmall } ] ) 
		\otimes \bigl( \sigma^{+}( \tau_1^{\ast} ) \boxtimes \sigma^{-}( \tau_2 ) \bigr) . 
\end{align*}
The theorem follows.
\end{proof}

Let us briefly describe our orbit correspondence in terms of signed Young diagrams.  
See also \cite{Ohta.preprint1}, \cite{Ohta.preprint2}, \cite{DKP.preprint}.

Since our $ \Gbig $ or $ \Gsmall $ is either $ Sp(2 n, \R) , O(p, q), U(m, n), Sp(p, q) $ or $ O^{\ast}(2 n) $, 
their nilpotent $ K_{\C} $-orbits are classified by the signed Young diagrams.  
See, e.g., \cite[\S 9.3]{Collingwood.McGovern.1993} for the classification in the case of classical Lie algebras.  
We reproduce it for the readers' convenience.

\begin{lemma}
The nilpotent $ K_{\C} $-orbits for a symmetric pair $ ( G, K ) $ are in one-to-one correspondence with the following signed Young diagrams.
\begin{thmenumerate}
\item
$ G = Sp(2 n, \R) $ : diagrams of size $ 2 n $ with any signature in which odd rows begin with $ + $ and occur with even multiplicity.  
\item
$ G = O(p, q) $ : diagrams of size $ p + q $ with signature $ ( p, q ) $ 
in which even rows begin with $ + $ and occur with even multiplicity.  
\item
$ G = U( m, n ) $ : diagrams of size $ m + n $ with signature $ ( m, n ) $.  
\item
$ G = O^{\ast}( 2 n ) $ : diagrams of size $ n $ with any signature in which odd rows begin with $ + $.
\item
$ G = Sp(p, q) $ : diagrams of size $ p + q $ with signature $ ( p, q ) $ in which even rows begin with $ + $.
\end{thmenumerate}
\end{lemma}

\begin{proposition}
\label{prop:signed.Y.diagram}
Let $ \Osmall $ be a nilpotent $ \KsmallC $-orbit in $ \nilpotents( \ssmall ) $, and 
$ \Obig = \theta( \Osmall ) $ its theta lift.  
If $ \Osmall $ corresponds to a signed Young diagram $ \lambda' $, 
$ \Obig $ corresponds to $ \lambda $ which is obtained by adding one box to the right end of each row (including empty one) of $ \lambda' $.  
The signature of the added boxes are completely determined by the requirement that 
$ \lambda $ is a signed Young diagram for $ ( \Gbig , \Kbig ) $.
\end{proposition}

\begin{proof}
We note that the stable range condition ensures that we can always make $ \lambda $ out of $ \lambda' $ by adding a box in each row.  
Then the proof is based on a case-by-case calculation, using explicit realization of moment maps in the Appendix.  
Since we know that $ \varphi( \psi^{-1}( \closure{\Osmall} ) ) = \closure{\Obig} $, we only need to take a representative of $ \Osmall $ and 
examine a generic element in $ \varphi( \psi^{-1}( \closure{\Osmall} ) ) $.  
\end{proof}

\section{Lifting of holomorphic nilpotent orbits}
\label{section:lift.hol.nilpotent.orbit}

We apply \theoremref{thm:fun.ring.of.lifted.orbit} to the lifting of the so-called holomorphic nilpotent orbits.  
In the process, we get a family of spherical orbits and prove normality of the closure of certain lifted orbits.  
This application reproduces some of our previous results 
(\cite{Nishiyama.MA.2000}, \cite{Nishiyama.JSPS.DFG.proc}, \cite{Nishiyama.Zhu.2001}) as well.  
In the following, we use the notations in \theoremref{thm:fun.ring.of.lifted.orbit} freely.

Let us consider the orbit decomposition of $ \ssmallp $ by the adjoint action of $ \KsmallC $.  
Any $ \KsmallC $-orbit in $ \ssmallp $ is clearly nilpotent.  
We call these nilpotent orbits \textit{holomorphic}.  
It is well known that $ \ssmallp $ is a prehomogeneous vector space, and 
there exists a numbering of $ \KsmallC $-orbits 
$ \Osmall_0 , \Osmall_1, \ldots, \Osmall_l $ in such a way that 
$ \Osmall_{i - 1} \subset \closure{ \Osmall_{i} } $ for $ 1 \leq i \leq l $.  
Here $l$ is the real rank of $\Gsmall $. 
As a consequence, 
$ \Osmall_0 =  \{ 0 \} $ and $ \closure{ \Osmall_l } = \ssmallp $, i.e., 
$ \Osmall_l $ is the open dense orbit in $ \ssmallp $.  
The orbit $ \Osmall_l $ is called \textit{regular}, while the orbits
$ \{ \Osmall_k \}_{0 \leq k < l}$ are called \textit{singular}.

\subsection{Trivial orbit}

We first consider the trivial orbit $ \Osmalltrivial = \Osmall_0 =\{ 0 \} $.  
Then, the lifted orbit $ \Obigtrivial = \theta( \Osmalltrivial ) $ is a two-step nilpotent orbit in 
$ \nilpotents( \sbig ) $ (cf.~\cite{Nishiyama.MA.2000}).  
The orbit $ \Obigtrivial $ corresponds to the following signed Young diagrams.
\begin{equation*}
\begin{array}{l@{\; : \quad }l}
O(p, q) & \lambda = [ ({+}{-})^n ({-}{+})^n ({+})^{p - 2 n} ({-})^{q - 2 n} ] \\
U(p, q) & \lambda = [ ({+}{-})^m ({-}{+})^n ({+})^{p - (m + n)} ({-})^{q - (m + n)} ] \\
Sp(p, q) & \lambda = [ ({+}{-})^n ({+})^{p - n} ({-})^{q - n} ] 
\end{array}
\end{equation*}
We have

\begin{theorem}
The group $ \KbigC $ acts on $ \closure{ \Obigtrivial } $ multiplicity-freely. As a $ \KbigC = \KbigpC \times \KbigqC $-module, we have
\begin{equation}
\C[ \closure{ \Obigtrivial } ] \simeq 
\directsum_{ \tau \in R( \KsmallC ) } \sigma^{+}( \tau^{\ast} ) \boxtimes \sigma^{-}( \tau ) .
\end{equation}
\end{theorem} 

\begin{proof}
Since $ \C[ \Osmalltrivial ] = \C $, by \theoremref{thm:fun.ring.of.lifted.orbit}, 
we get 
\begin{displaymath}
\C[ \closure{ \Obigtrivial } ] \simeq 
\directsum_{ \tau_1, \tau_2 \in R( \KsmallC ) } \Hom_{\KsmallC}( \tau_1 \otimes \tau_2^{\ast}, \C ) \otimes 
\bigl( \sigma^{+}( \tau_1^{\ast} ) \boxtimes \sigma^{-}( \tau_2  ) \bigr) .
\end{displaymath}
But, by Schur's lemma, the multiplicity $ \Hom_{\KsmallC}( \tau_1 \otimes \tau_2^{\ast}, \C ) $ is not zero if and only if $ \tau_1 = \tau_2 $, in which case it is $ \C $, i.e., the multiplicity is one.
\end{proof}

Recall that an orbit $ \Obig $ is called \textit{spherical}, or more precisely  
\textit{$ \KbigC $-spherical}, if a Borel subgroup of $ \KbigC $ 
has a dense orbit in $ \Obig $. This is equivalent to saying that 
$ \KbigC $ acts on $ \closure{ \Obig } $ multiplicity-freely (i.e., 
$ \C[ \closure{ \Obig } ] $ decomposes without multiplicity as a $ \KbigC $-module).  

\begin{corollary}
The orbit $ \Obigtrivial $ lifted from the trivial orbit is spherical.  
The closure $ \closure{\Obigtrivial} $ is a normal variety.
\end{corollary}

\begin{proof}
The claim that $ \Obigtrivial $ is spherical follows directly from the above theorem.  
The normality follows from \cite[Th.~10]{Vinberg.Popov.1972}.
\end{proof}

\subsection{Regular holomorphic orbit}

Let us denote by $ \Obighol_l = \theta( \Osmall_l ) $ the theta lift from the regular holomorphic orbit $ \Osmall_l \subset \ssmallp $.  
The orbit $ \Obighol_l $ corresponds to the following signed Young diagram.  
\begin{equation*}
\begin{array}{l@{\; : \quad }l}
O(p, q) & \lambda = [ ({+}{-}{+})^n ({+})^{p - 2 n} ({-})^{q - n} ] \\
U(p, q) & \lambda = [ ({+}{-}{+})^n ({+}{-})^{m - n} ({+})^{p - (m + n)} ({-})^{q - m} ] \quad ( m \geq n ) \\[3pt]
Sp(p, q) & \lambda = 
\left\{ \begin{array}{l}
{[ ({+}{-}{+})^{(n - 1)/2} ({+}{-}) ({+})^{p - n} ({-})^{q - (n + 1)/2} ]} \quad ( n = 2 l + 1 ) \\
{[ ({+}{-}{+})^{n/2} ({+})^{p - n} ({-})^{q - n/2} ]} \quad ( n = 2 l ) 
\end{array} \right. 
\end{array}
\end{equation*}
For the regular orbit, we have 
\begin{displaymath}
\Xi( \Osmall_l ) = \momentsmall^{-1}( \closure{ \Osmall_l } ) 
	= \momentsmall^{-1}( \ssmallp ) = \Wp \times \nullconeq, 
\end{displaymath}
where $\momentsmallq =\momentsmall \restrict _{\Wq}$ and 
$ \nullconeq = ( \momentsmallq )^{-1}( 0 ) \subset \Wq $ is a null cone.  
Thus \theoremref{thm:fun.ring.of.lifted.orbit} yields the following

\begin{lemma}
Let $ \Obighol_l \subset \nilpotents( \sbig ) $ be the theta lift of the regular holomorphic orbit $ \Osmall_l $, 
which is open dense in $ \ssmallp $.  
Then its closure is an affine quotient  of $ \Wp \times \nullconeq $ by $ \KsmallC $, 
\begin{displaymath}
\closure{\Obighol_l} \simeq ( \Wp \times \nullconeq ) \GITquotient \KsmallC .
\end{displaymath}
\end{lemma}

\mycomment{
On the claim that $ (\text{vector space}) \times (\text{normal variety}) $ is again normal: 

This can be proved as follows (I have learned this from Makoto Matumoto).  
What should be proved is the following algebraic statement.  

\begin{myproposition}
Let $ A $ be an integrally closed domain over $ \C $.  
Then the polynomial algebra $ A [ x_1 , \ldots , x_n ] $ over $ A $ is integrally closed.  
\end{myproposition}

To prove this, it is enough to assume that $ n = 1 $.  
The following lemmas are well known.

\begin{mylemma}
If $ A $ is a discrete valuation ring, then $ A[ x ] $ is a UFD.  
In particular, $ A[ x ] $ is integrally closed.
\end{mylemma}

\begin{mylemma}
If $ A $ is an integrally closed domain, then 
\begin{math}
A = \bigcap_{\lie{p} : \text{ \upshape height 1 }} A_{\lie{p}} .
\end{math}
\end{mylemma}

Let $ K = Q( A ) $ be the quotient field of $ A $.  
By the first lemma, we know that $ A_{\lie{p}} [ x ] $ is integrally closed in $ Q( A_{\lie{p}} ) ( x ) $.  
Since the integral closure of $ A [ x ] $ in $ K( x ) $ is contained in the integral closure of 
$ A_{\lie{p}} [ x ] $ in $ K( x ) = Q( A_{\lie{p}} ) ( x ) $ for any $ \lie{p} $, we have
\begin{displaymath}
\text{int.~closure}( A [ x ] ) \subset \bigcap_{ \lie{p} : \text{ height 1 } } A_{\lie{p}} [ x ] = A [ x ] .
\end{displaymath}
This completes the proof of the proposition.
}

Let us recall ($\S $\ref{sec:diamond}) 
the compact group $ \Lbigp $ whose action on $ \Wp $ commutes with that of $ \Ksmall $.  
In fact $ ( \Lbigp , \Ksmall ) $ forms a dual pair in $ Sp( \Wp_{\R} ) $. 
Let  
\begin{equation}
\label{eqn:Lbigp.Ksmall.duality.original}
\C[ \Wp ] \simeq \directsum_{ \tau \in \irreps{ \KsmallC ; \C[ \Wp ] } 
	} 
	\rho^{+}( \tau ) \otimes \tau 
\end{equation}
be the decomposition of $ \C[ \Wp ]$ as an $ \Lbigp \times \Ksmall $-module.  
Here $ \rho^{+}( \tau ) \in \irreps{ \Lbigp } $ corresponds to  
$ \tau \in \irreps{ \Ksmall } $ via the above multiplicity-free decomposition. 
Note that $\C[ \Wp ] \simeq \harmonics( \KbigpC ) \otimes \C[ \Wp ]^{\KbigpC}$,
and so $ \irreps{ \KsmallC ; \C[ \Wp ] } =
\irreps{ \KsmallC ; \harmonics(\KbigpC) } = R( \KsmallC ) ^{\ast} $.
We may therefore rewrite Equation (\ref{eqn:Lbigp.Ksmall.duality.original}) as
\begin{equation}
\label{eqn:Lbigp.Ksmall.duality}
\C[ \Wp ] \simeq \directsum_{ \tau \in R( \KsmallC ) } 
	\rho^{+}( \tau ^{\ast}) \otimes \tau ^{\ast}. 
\end{equation}

Note that $ \Lbigp $ is a unitary group containing $ \Kbigp $ in each of the three cases, 
hence its complexification $ \LbigpC $ is isomorphic to a general linear group.  
Since the action of $ \LbigpC $ commutes with that of $ \KsmallC $, $ \LbigpC \times \KbigqC $ naturally acts on the affine quotient space 
\begin{displaymath}
( \Wp \times \nullconeq ) \GITquotient \KsmallC \simeq \closure{ \Obighol_l }, 
\end{displaymath}
extending the action of $ \KbigC = \KbigpC \times \KbigqC $.  
Note that the orbit $ \Obighol_l $ itself does \textit{not} admit an action of $ \LbigpC \times \KbigqC $, but its closure does.
 
\begin{theorem} 
The group $ \LbigpC \times \KbigqC $ acts on the closure $ \closure{ \Obighol_l } $ multiplicity-freely. 
As a $ \LbigpC \times \KbigqC $-module, we have
\begin{equation}
\C[ \closure{ \Obighol_l } ] 
\simeq \directsum_{ \tau \in R( \KsmallC ) } \rho^{+}( \tau^{\ast} ) \boxtimes \sigma^{-}( \tau ).
\label{eq:sph.decomp.closure.Obighol}
\end{equation}
Consequently, the closure $ \closure{\Obighol_l} $ is $ \LbigpC \times \KbigqC $-spherical, and it is a normal variety.
\end{theorem}

\begin{proof} We have 
\begin{align}
\C[ \closure{ \Obighol_l } ] 
&\simeq \bigl( \C[ \Wp ] \otimes \C[ \nullconeq ] \bigr)^{\KsmallC} \notag \\
&\simeq \Bigl( \Bigl( \directsum_{ \tau _1\in R( \KsmallC ) } 
	\rho^{+}( \tau _1^{\ast}) \otimes \tau _1^{\ast} \Bigr) \otimes 
	 \Bigl( \directsum_{ \tau_2 \in R( \KsmallC ) } \sigma^{-}( \tau_2 ) \otimes \tau_2 \Bigr) \Bigr)^{\KsmallC}  \notag \\
&\simeq \directsum_{ \tau \in R( \KsmallC ) } \rho^{+}( \tau^{\ast} ) \boxtimes \sigma^{-}( \tau ).\notag
\end{align}
The normality follows from \cite[Th.~10]{Vinberg.Popov.1972}, 
since $ \Z R( \KsmallC ) \cap \Q_+ R( \KsmallC ) = R( \KsmallC ) $.
\end{proof}

\subsection{Degree of nilpotent orbits lifted from holomorphic orbits}

As an application of the explicit decomposition formula of the function ring of nilpotent orbits, 
we give an integral formula of their degrees.   

Since a nilpotent orbit $ \orbit \subset \sbig $ is a cone, we naturally consider it in a projective space $ \bbP( \sbig ) $.  As a projective variety, the closure $ \bbP( \closure{\orbit} ) $ in $ \bbP( \sbig ) $ has a degree. We shall denote this degree by $ \deg \closure{\orbit} $ by a slight abuse of notation, and we call it the degree of the nilpotent orbit $ \orbit $. We refer to \cite[\S V.4.2]{Shafarevich.1994} for the definition of the degree of projective varieties. See also \cite{Hartshorne.1977}.  

To state our result, let us introduce some notation.  
We denote by $ \Delta = \Delta(X_r) $ an irreducible root system of type $ X_r = A_{r - 1}, B_r, C_r, D_r $.  
We realize $ \Delta $ in a standard way (see \cite[Planche I-IX]{Bourbaki.1968}), and choose a positive system $ \Delta^+ $ as 
$ \Delta^+(A_{r - 1}) = \{ \varepsilon_i - \varepsilon_j \mid 1 \leq i < j \leq r \} $; 
$ \Delta^+(D_r) = \{ \varepsilon_i \pm \varepsilon_j \mid 1 \leq i < j \leq r \} $; 
$ \Delta^+(B_r) = \Delta^+(D_r) \cup \{ \varepsilon_i \mid 1 \leq i \leq r \} $; 
$ \Delta^+(C_r) = \Delta^+(D_r) \cup \{ 2 \varepsilon_i \mid 1 \leq i \leq r \} $.  
As usual, $ \rho = \rho_{X_r} = \sum_{\alpha \in \Delta^+} \alpha / 2 $ denotes the half sum of positive roots, and   
we fix an invariant inner product $ \langle \alpha , \beta \rangle $ on the root space 
which makes $ \{ \varepsilon_i \} $ an orthonormal system.  
For fixed $ n \leq r $, we put 
\begin{align*}
\Delta^+_n( X_r ) &= \{ \alpha \in \Delta^+(X_r) \mid \langle \alpha , \varepsilon_i \rangle = 0 \; ( 1 \leq i \leq n ) \} , 
\; \text{ and } \; \\
\Phi^+_n(X_r) &= \Delta^+( X_r ) \setminus \Delta^+_n( X_r ) ,
\end{align*}
and for $ m+n \leq r$, 
\begin{displaymath}
\Phi^+_{m, n}(A_{r - 1}) = \{ \varepsilon_i - \varepsilon_j \in \Delta^+(A_{r - 1}) \mid i \leq m \text{ or } r - n < j \} .
\end{displaymath}

We introduce the difference product $ D_n(x) $ on a multi-variable $ x = ( x_{1}, \ldots, x_{n} ) $ defined by 
\begin{equation*}
\textstyle D_n(x) = \prod\nolimits_{1 \leq i < j \leq n} ( x_i - x_j ) .
\end{equation*}
We abbreviate $ D_n(x^2) = D_n( x_1^2, \ldots, x_n^2 ) $. Finally set
\begin{displaymath}
\textstyle \Omega_n = \{ x = ( x_i )_{1 \leq i \leq n} \mid x_i \geq 0 , \sum\nolimits_{i = 1}^n x_i \leq 1 \}.
\end{displaymath}

\begin{proposition}
\label{prop:degree.of.orbits.integral.expression}
\begin{thmenumerate}
\item
Case $ O(p, q) \times Sp(2 n, \R) \; ( 2 n < p , q ) $: 
Put $ r = [ p / 2 ], s = [ q / 2 ] $ and $ X = B $ or $ D $ (respectively $ Y = B $ or $ D $) according as $ p $ (respectively $ q $) 
is odd or even.  
Then the degree of the nilpotent orbit is given by 
\begin{align*}
\deg \closure{\Obigtrivial} &= \frac{1}{n!} {\textstyle
\prod\limits_{\alpha \in \Phi^+_n(X_r)} \langle \rho_{X_r}, \alpha \rangle^{-1}
\prod\limits_{\beta \in \Phi^+_n(Y_s)} \langle \rho_{Y_s}, \beta \rangle^{-1} }
\int_{\Omega_n} |D_n(x^2)|^2 \, { \textstyle \prod\limits_{i = 1}^n x_i^{p + q - 2 n} } dx  \\
\deg \closure{\Obighol_l} &= \frac{1}{n!}  { \textstyle
\prod\limits_{\alpha \in \Phi^+_n(X_r)} \langle \rho_{X_r}, \alpha \rangle^{-1}
\prod\limits_{\beta \in \Phi^+_n(A_{q - 1})} \langle \rho_{A_{q - 1}}, \beta \rangle^{-1} }
\int_{\Omega_n} |D_n(x^2) D_n(x)| \, { \textstyle \prod\limits_{i = 1}^n x_i^{p + q - 3 n } } dx  
\end{align*}
\item
Case $ U(p, q) \times U(m, n) \; ( m + n \leq p, q ) $:
{\small
\begin{align*}
\deg \closure{\Obigtrivial} &= \frac{1}{m! n!} \textstyle
\prod\limits_{\alpha \in \Phi^+_{m, n}(A_{p - 1})} \langle \rho_{A_{p - 1}}, \alpha \rangle^{-1}
\prod\limits_{\beta \in \Phi^+_{m, n}(A_{q - 1})} \langle \rho_{A_{q - 1}}, \beta \rangle^{-1} \\
& \times 
\int\limits_{\Omega_m \times \Omega_n}  
|D_m(x)|^2 |D_n(y)|^2 
\textstyle
\prod\limits_{\begin{smallmatrix}
1 \leq i \leq m\\
1 \leq j \leq n
\end{smallmatrix}} ( x_i + y_j )^2 
\prod\limits_{i = 1}^m x_i^{ p + q - 2 ( m + n ) } \prod\limits_{j = 1}^n y_j^{ p + q - 2 ( m + n ) } dx dy 
\end{align*}
\begin{align*}
\deg \closure{\Obighol_l} &= \frac{1}{m! n!} \textstyle
\prod\limits_{\alpha \in \Phi^+_{m, n}(A_{p - 1})} \langle \rho_{A_{p - 1}}, \alpha \rangle^{-1} 
\prod\limits_{\beta \in \Phi^+_{m}(A_{q - 1})} \langle \rho_{A_{q - 1}}, \beta \rangle^{-1} 
\prod\limits_{\gamma \in \Phi^+_{n}(A_{q - 1})} \langle \rho_{A_{q - 1}}, \gamma \rangle^{-1} \\
&\times 
\int\limits_{\Omega_m \times \Omega_n} 
|D_m(x)|^2 |D_n(y)|^2 
\textstyle
\prod\limits_{\begin{smallmatrix}
1 \leq i \leq m\\
1 \leq j \leq n
\end{smallmatrix}} ( x_i + y_j )  
\prod\limits_{i = 1}^m x_i^{ p + q - ( 2 m + n ) } \prod\limits_{j = 1}^n y_j^{ p + q - ( m + 2 n ) } dx dy 
\end{align*}
}
\item
Case $ Sp(p, q) \times O^{\ast}(2n) \  ( n \leq p , q ) $: 
\begin{align*}
\deg \closure{\Obigtrivial} &= \frac{2^{2 n}}{n!} { \textstyle
\prod\limits_{\alpha \in \Phi^+_n(C_p)} \langle \rho_{C_p}, \alpha \rangle^{-1}
\prod\limits_{\beta \in \Phi^+_n(C_q)} \langle \rho_{C_q}, \beta \rangle^{-1} } 
\int_{\Omega_n} |D_n(x^2)|^2 \, { \textstyle \prod\limits_{i = 1}^n x_i^{2 ( p + q ) - 4 n + 2 } } dx  \\
\deg \closure{\Obighol_l} &= \frac{2^{n}}{n!}  { \textstyle
\prod\limits_{\alpha \in \Phi^+_n(C_p)} \langle \rho_{C_p}, \alpha \rangle^{-1}
\prod\limits_{\beta \in \Phi^+_n(A_{2 q - 1})} \langle \rho_{A_{2 q - 1}}, \beta \rangle^{-1} } \\
& \hspace*{.25\textwidth} \times 
\int_{\Omega_n} |D_n(x^2) D_n(x)| \, { \textstyle \prod\limits_{i = 1}^n x_i^{2 ( p + q ) - 3 n + 1 } } dx  
\end{align*}
\end{thmenumerate}
\end{proposition}

\begin{proof}
Since we know the $ \KbigC $-decomposition of the regular function ring, 
we can formally express the Poincar\'{e} series by using Weyl's dimension formula.  
The degree $ e $ is encoded in the top term of the Hilbert polynomial $ H_{\closure{\Obig}}(t) $, 
whose degree coincides with $ d = \dim \bbP( \closure{\Obig} ) $.  
Thus, we have $ e = d! \lim_{t \rightarrow \infty} H_{\closure{\Obig}}(t) / t^d $, which is essentially expressed by a Riemann integral (the limit of a Riemann sum).  
This machinery is explained in \cite[\S 2.2]{KO.2001}. We omit the details.
\end{proof}

\newcommand{\sign}{\mbox{sgn}}
\newcommand{\vdM}{\mbox{vdM}}

Let us evaluate the following integral, which appears in the above proposition.

\begin{theorem}
Let $\kappa$ be a complex number with a positive real part.
Then
\begin{equation}
\frac{1}{n!} \int_{\Omega_n} D_n(x^2) D_n(x) \prod_{i=1}^n x_i^{\kappa-1} dx
= \frac{2^{n(n-1)/2} \prod_{i=0}^{n-1} i! \, \Gamma(\kappa+2i)}%
{\Gamma(3n(n-1)/2 + n\kappa +1)}.
\label{eq:DD}
\end{equation}
\end{theorem}

\begin{proof} First we set
\[
\textstyle f(x)=f(x,\kappa) = D_n(x^2) D_n(x) \prod\nolimits_{i=1}^n x_i^{\kappa-1} .
\]
Let $a=3n(n-1)/2 + n\kappa$, then $f(x)$ is homogeneous of degree $a-n$,
that is,
\[
f(\lambda x) = \lambda^{a-n} f(x)
\qquad\quad (\lambda >0).
\]
We consider the integral
\[
\Gamma(a+1) \int_{\Omega_n} f(x) dx
= \int_{\Omega_n\times (0,\infty)} f(x) s^a e^{-s} dx ds.
\]
We have a diffeomorphism
\[
\Omega_n \times (0,\infty) \ni (x,s)
\mapsto
(y,t) \in (0,\infty)^n \times (0,\infty),
\]
given by
\[
\left\{
\begin{array}{ccl}
x_i &=& y_i/(t+\sum_{k=1}^n y_k), \\
s &=& t+ \sum_{k=1}^n y_k,
\end{array}
\right.
\qquad
\left\{
\begin{array}{ccl}
y_i &=& s x_i, \\
t &=& (1-\sum_{i=1}^n x_i)s.
\end{array}
\right.
\]
The Jacobian of this map is given by
\[
\frac{\partial(y,t)}{\partial(x,s)} = s^n.
\]
Then
\[
f(x) s^a e^{-s} dx ds
= f(sx) s^n e^{-s} dx ds
= f(y) e^{-t-\sum_{i=1}^n y_i} dy dt.
\]
This proves
\begin{eqnarray}
\Gamma(a+1) \int_{\Omega_n} f(x) dx
&=& \int_{\Omega_n\times(0,\infty)} f(x) s^a e^{-s} dx ds  = \int_{(0,\infty)^{n+1}} f(y) e^{-\sum_{i=1}^n y_i} e^{-t}  dy dt \notag \\
&=& \int_{(0,\infty)^{n}} f(y) e^{-\sum_{i=1}^n y_i} dy
\label{shiki:limit}
\end{eqnarray}
Now we use the determinantal expression
\begin{eqnarray*}
D_{n}(y^2) &=& \prod_{1\le i<j\le n}(y_i^2-y_j^2)
= \det(y_j^{2(i-1)})_{i,j=1,\dots,n}
= \sum_{\sigma \in S_n} \sign(\sigma) \prod_{i=1}^n y_i^{2(\sigma(i)-1)}, \\
D_n(y) &=& \prod_{1\le i<j\le n}(y_i-y_j)
= \det(y_j^{i-1})_{i,j=1,\dots,n}
= \sum_{\tau \in S_n} \sign(\tau) \prod_{i=1}^n y_i^{\tau(i)-1}.
\end{eqnarray*}
Then (\ref{shiki:limit}) becomes 
\begin{align*}
\int_{(0,\infty)^{n}} f(y) & e^{-\sum_{i=1}^n y_n} dy
=
\sum_{\sigma,\tau} \sign(\sigma)\sign(\tau)
\prod_{i=1}^n \int_0^\infty
y_i^{2(\sigma(i)-1)+(\tau(i)-1)+\kappa-1}
e^{- y_i} dy_i
\\
&=
\textstyle 
\sum\nolimits_{\sigma,\tau}  \sign(\sigma) \, \sign(\tau)
\prod\nolimits_{i=1}^n
\Gamma(2\sigma(i)+\tau(i)+\kappa-3) \\
&=
\textstyle 
n!
\sum\nolimits_{\sigma}  \sign(\sigma)
\prod\nolimits_{i=1}^n
\Gamma(2\sigma(i)+i+\kappa-3)  \\
&=
\textstyle 
n!
\det(\Gamma(2i+j+\kappa-3))_{i,j=1,\dots,n} \\[5pt]
&=
\textstyle 
n! \prod\nolimits_{i=1}^n \Gamma(2i+\kappa-2)\times
\displaystyle
\det\left(\frac{\Gamma(2i+j+\kappa-3)}{\Gamma(2i+\kappa-2)}\right)_{i,j=1,\dots,n}
\\
&=
\textstyle 
n! \prod\nolimits_{i=0}^{n-1} \Gamma(\kappa+2i)\times
\det\left((2i+\kappa-2)_{j-1} \right)_{i,j=1,\dots,n},
\end{align*}
where $(c)_k=c(c+1)\cdots (c+k-1)$. The last determinant is
\begin{align*}
\det &\left( (2i+\kappa-2)_{j-1} \right)_{i,j=1,\dots,n}
=
\det\left((2i+\kappa-2)^{j-1} \right)_{i,j=1,\dots,n}
\\
&=
\vdM(\kappa,\kappa+2,\dots,\kappa+2n-2) 
= \vdM(0,2,\dots,2n-2) 
= 2^{n(n-1)/2} 
\textstyle 
\prod_{i=1}^{n-1} i!
\end{align*}
Here $\vdM$ denotes the vandermonde determinant:
\[
\vdM(c_1,\dots,c_n)
= \textstyle \prod_{1\le i < j \le n} (c_j-c_i).
\]
This gives 
\[
\frac{1}{n!} \int_{\Omega_n} f(x) dx
= \frac{1}{\Gamma(a+1)} 2^{n(n-1)/2} \prod_{i=0}^{n-1} i! \, \Gamma(\kappa+2i) , 
\]
which is the formula (\ref{eq:DD}). 
\end{proof}

\begin{remark} By a similar argument, we have
\[
\frac{1}{n!} \int_{\Omega_n} D_n(x^2)^2 \prod_{i=1}^n x_i^{\kappa-1} dx
= \frac{\prod_{i=0}^{n-1} \Gamma(\kappa+2i)}%
{\Gamma(2n(n-1) + n\kappa +1)}
\times
\det\left((\kappa+2i-2)_{2j-2} \right)_{i,j=1,\dots,n}
\label{eq:D^2}
\]
The last determinant depends on $ \kappa $ polynomially and seems to have no simple expression.
This suggests that $\closure{\Obigtrivial}$ may have a more complicated structure
than that of $\closure{\Obighol_l}$.
\end{remark}

\section{Appendix}

In Appendix, we will give an explicit construction of moment maps $ \momentbig $ and $ \momentsmall $, and establish some basic properties of null cones and harmonics.  

\subsection{$ O(p, q) \times Sp(2 n, \R) \  ( 2 n < \mathoperator{min} ( p, q ) ) $}

This case is already treated in \cite{Nishiyama.MA.2000}.  
For convenience of the readers, we reproduce the synopsis of the arguments, and at the same time we improve some of statements there.

Put 
$ \Gbig = O(p, q) $ and $ \Gsmall = Sp(2 n, \R) $, 
and assume the stable range condition $ 2 n < p, q $.  We denote:
{\small
\begin{gather*}
\begin{array}{ccc}
\Mbig = U(p, q) & \supset & \Gbig = O(p, q) \\
\cup & & \cup \\
\begin{array}[t]{l@{\,}l}
\Lbig &= \Lbigp \times \Lbigq \\
      &= U(p) \times U(q) 
\end{array}
& \supset & 
\begin{array}[t]{l@{\,}l}
\Kbig &= \Kbigp \times \Kbigq \\
      &= O(p) \times O(q)
\end{array}
\end{array}
\qquad 
\begin{array}{ccc}
\Msmall = Sp(2 n, \R)^2 & \overset{\Delta}{\supset} & \Gsmall = Sp(2 n, \R) \\
\cup & & \cup \\
\begin{array}[t]{l@{\,}l}
\Lsmall &= \Ksmall \times \Ksmall \\
        &= U(n)^2 
\end{array}
& \overset{\Delta}{\supset} & \Ksmall = U(n)
\end{array}
\end{gather*}
}
In this diagram, the vertical containment means respectively maximal compact subgroups, and $ \Delta $ denotes the diagonal embedding.  Put
\begin{align*}
W 
& = M_{p + q, n}(\C) = \left\{ \begin{pmatrix} A \\ B \end{pmatrix} \mid A \in M_{p, n}(\C) , B \in M_{q, n}(\C) \right\} \\
& = \Wp \oplus \Wq = M_{p, n}(\C) \oplus M_{q, n}(\C) .
\end{align*}
The complexification $ \KbigC = O(p, \C) \times O(q, \C) $ acts on $ W $ via the left multiplication, and 
$ \KsmallC = GL_n(\C) $ acts on $ W $ as 
\begin{displaymath}
g \cdot \begin{pmatrix} A \\ B \end{pmatrix} = \begin{pmatrix} A \transpose{g} \\ B g^{-1} \end{pmatrix} , \qquad
g \in GL_n(\C) .
\end{displaymath}
We fix a Cartan decomposition $ \gbig = \kbig \oplus \sbig $ (resp. $ \gsmall = \ksmall \oplus \ssmall $) as 
\begin{align*}
\gbig 
&= \mathfrak{o}( p + q , \C ) = 
	\left( \begin{array}{@{}c|c@{}}
	\rule[-1.3ex]{0pt}{3.3ex} \Alt_p(\C) & 0          \\ \hline
	\rule[-1ex]{0pt}{3.5ex}   0          & \Alt_q(\C) 
	\end{array} \right) 
\oplus 
	\left( \begin{array}{@{}c|c@{}}
	\rule[-1.3ex]{0pt}{3.3ex} 0                         & M_{p, q}(\C) \\ \hline
	\rule[-1ex]{0pt}{3.5ex}   \transpose{M}_{p, q}(\C)  & 0            
	\end{array} \right)
= \kbig \oplus \sbig , \\[10pt]
\gsmall 
&= \mathfrak{sp}( 2 n , \C ) = 
	\left( \begin{array}{@{}c|c@{}}
	\rule[-1.3ex]{0pt}{3.3ex} M_n(\C) & 0          \\ \hline
	\rule[-1ex]{0pt}{3.5ex}   0          & - \transpose{M}_n(\C) 
	\end{array} \right) 
\oplus 
	\left( \begin{array}{@{}c|c@{}}
	\rule[-1.3ex]{0pt}{3.3ex} 0           & \Sym_n(\C) \\ \hline
	\rule[-1ex]{0pt}{3.5ex}   \Sym_n(\C)  & 0            
	\end{array} \right)
= \ksmall \oplus \ssmall .
\end{align*}
Thus, we can identify 
$ \sbig = M_{p, q}(\C) $ and $ \ssmall = \ssmallp \oplus \ssmallq = \Sym_n(\C) \oplus \Sym_n(\C) $.  
The moment maps are explicitly given by 
\begin{displaymath}
\momentbig : W \ni \begin{pmatrix} A \\ B \end{pmatrix} \longmapsto A \transpose{B} \in \sbig , \qquad 
\momentsmall : W \ni \begin{pmatrix} A \\ B \end{pmatrix} \longmapsto ( \transpose{A} A , \transpose{B} B ) \in \ssmall .
\end{displaymath}
These maps are clearly $ \KbigC \times \KsmallC $-equivariant if we define the trivial $ \KbigC $-action (resp. $ \KsmallC $-action) on $ \ssmall $ (resp. $ \sbig $).  
In the language of algebra, $ \momentbig $ is induced by 
\begin{displaymath}
\momentbig^{\ast}( z_{i, j} ) = \sum_{k = 1}^n a_{i, k } b_{j, k} , \qquad 
	Z = ( z_{i, j} ) \in \sbig , \quad 
	\begin{pmatrix} A \\ B \end{pmatrix} = \begin{pmatrix} ( a_{i, j} ) \\ ( b_{i, j} ) \end{pmatrix} \in W , 
\end{displaymath}
and it is well known that $ \momentbig^{\ast}( z_{i, j} ) $'s generate $ GL_n $-invariants $ \C[ W ]^{\KsmallC} $.  
This shows that $ \momentbig : W \rightarrow \sbig $ is an affine quotient map onto its closed image:  
\begin{math}
\sbig \supset \momentbig( W ) \simeq W \GITquotient \KsmallC .  
\end{math}
Similarly, $ \momentsmall $ is induced by 
\begin{displaymath}
\textstyle
\momentsmall^{\ast}( x_{i, j} ) = \sum_{k = 1}^p a_{k, i} a_{k, j} , \; X = ( x_{i, j} ) \in \ssmallp \, ; \quad 
\momentsmall^{\ast}( y_{i, j} ) = \sum_{l = 1}^q b_{l, i} b_{l, j} , \; Y = ( y_{i, j} ) \in \ssmallq .
\end{displaymath}
Again, $ \momentsmall( x_{i, j} )^{\ast} $'s and $ \momentsmall( y_{i, j} )^{\ast} $'s generate the $ O(p, \C) \times O(q, \C) $-invariants 
$ \C[ W ]^{\KbigC} $, which proves that 
$ \momentsmall : W \rightarrow \ssmall $ is an affine quotient map.  
By the assumption of the stable range, $ 2 n < p, q $, it is easy to see that $ \momentsmall : W \rightarrow \ssmall $ is a surjection, hence
\begin{math}
\ssmall \simeq W \GITquotient \KbigC .
\end{math}

Let us consider the null cones
\begin{displaymath}
\nullcone = \momentsmall^{-1}( 0 ) \subset W , \qquad 
\nullconepq = ( \momentsmallpq )^{-1}( 0 ) \subset \Wpq , 
\end{displaymath}
where $ \momentsmallpq : \Wpq \rightarrow \ssmallpq $ is the restriction of $ \momentsmall $ to $ \Wpq $.

\begin{proposition}
\label{prop:OSp.nullcone}
Assume the stable range condition $ 2 n < p, q $.
\begin{thmenumerate}
\item
\label{item1:prop:OSp.nullcone}
The null cone $ \nullcone = \nullconep \times \nullconeq $ is an irreducible normal variety, which is of complete intersection 
in the set theoretic sense.  
\item
\label{item2:prop:OSp.nullcone}
There are finitely many $ \KbigC \times \KsmallC $-orbits in $ \nullcone $, which are completely classified by the ranks of each component in $ \nullconepq $.  
Among them, there is an open dense $ \KbigC \times \KsmallC $-orbit $ \sorbit_{n,n} $, which is a single $ \KbigC $-orbit.    
Moreover, the singular locus coincides with $ \nullcone \backslash \sorbit_{n, n} $, which is of codimension $ \geq 2 $.
\item
\label{item3:prop:OSp.nullcone}
The regular function ring $ \C[ \nullcone ] $ is isomorphic to the harmonics $ \harmonics( \KbigC ) \simeq \harmonicsp( \KbigpC ) \otimes \harmonicsq( \KbigqC ) $ as a $ \KbigC \times \KsmallC $-module.  
Its decomposition is explicitly given as
\begin{displaymath}
\C[ \nullconep ] \simeq \harmonicsp( \KbigpC ) 
	\simeq \directsum\limits_{\lambda \in \partition_n} \sigma_{\lambda}^{(p)}{}^{\ast} \otimes \tau_{\lambda}^{(n)}{}^{\ast} , 
\quad 
\C[ \nullconeq ] \simeq \harmonicsq( \KbigqC ) 
	\simeq \directsum\limits_{\lambda \in \partition_n} \sigma_{\lambda}^{(q)} \otimes \tau_{\lambda}^{(n)} , 
\end{displaymath}
where $ \partition_n $ denotes the set of all the partitions of length $ \leq n $, and 
$ \tau_{\lambda}^{(n)} \  (\text{resp. $ \sigma_{\lambda}^{(p)} $}) $
is an irreducible finite dimensional representation of 
$ GL_n \  (\text{resp. $ O(p, \C) $}) $ of highest weight $ \lambda $.
\end{thmenumerate}
\end{proposition}

\begin{remark}
\begin{thmenumerate}
\item
In \cite{Nishiyama.MA.2000}, the orbit decomposition of the null cone is carried out for $ \KbigC \times \LsmallC $, 
though it is more natural to consider $ \KbigC \times \KsmallC $-orbits.  
In fact, their orbits are the same (under the assumption of stable range).  
Also, in \cite{Nishiyama.MA.2000}, the normality is proved when $ 2n + 2 \leq p, q $.  
Here we improve it, and there is no more restriction other than the stable range condition.
\item
It is rather subtle to specify the representation $ \sigma_{\lambda}^{(p)} $ because $ O(p, \C) $ is not connected.  
For this, see \cite[\S 3.6.2]{Howe.1995.SchurLecture}.
\end{thmenumerate}
\end{remark}

\begin{proof}
It is straightforward to check the claim for the $ \KbigC \times \KsmallC $-orbit decomposition in \eqnref{item2:prop:OSp.nullcone} directly by calculation.  
Let us denote $ \sorbit_{r, s} = \sorbitp_r \times \sorbitq_s \in \nullcone $, where 
\begin{align*}
\sorbitp_{r} &= \{ A \in M_{p, n} \mid \transpose{A} A = 0 , \rank A = r \} \subset \nullconep , \\
\sorbitq_{s} &= \{ B \in M_{q, n} \mid \transpose{B} B = 0 , \rank B = s \} \subset \nullconeq .
\end{align*}
Then $ \sorbitp_r $ is a $ \KbigpC \times \KsmallC $-orbit, and $ \sorbit_{r,s} $ is a single $ \KbigC \times \KsmallC $-orbit.  
Since they are classified by rank, the orbit of the largest possible rank $ \sorbit_{n, n} $ is open dense in $ \nullcone $.  
This implies that $ \nullcone = \closure{ \sorbit_{n, n} } $ is irreducible.  
In fact, $ \sorbit_{n,n} $ is a single $ \KbigC $-orbit.  
This follows from the Witt theorem.  

Explicit calculation of the rank of differentials of the defining equations of 
$ \momentsmall^{\ast}( x_{ij} ) $'s and $ \momentsmall^{\ast}( y_{ij} ) $'s tells us that a point from $ \nullcone \backslash \sorbit_{n, n} $ is singular.  
Since 
\begin{displaymath}
\dim \nullconep = n p  - \frac{n (n + 1)}{2} , \qquad 
\dim \nullconeq = n q  - \frac{n (n + 1)}{2} , 
\end{displaymath}
we have $ \codim \nullcone = n ( n + 1 ) $, which is equal to the number of the defining equations.  
This proves that $ \nullcone $ is of complete intersection.  

The dimension of the $ \KbigpC \times \KsmallC $-orbit in $ \nullconep $ of rank $ r $ is given by 
\begin{displaymath}
\dim \sorbitp_r = r ( p + n ) - r^2 - \frac{ r ( r + 1 ) }{2} .
\end{displaymath}
From this formula, we can easily show that $ \codim ( \nullcone \backslash \sorbit_{n, n} ) \geq 2 $.  

By Kostant theorem (see \cite[Th.~2.2.11]{Chriss.Ginzburg.1997}, for example), 
we know that the defining radical ideal of $ \nullcone $ is generated by 
$ \momentsmall^{\ast}( x_{ij} ) $'s and $ \momentsmall^{\ast}( y_{ij} ) $'s, which are basic invariants of $ O(p, \C) \times O(q, \C) $.  This means that 
\begin{displaymath}
\C[ \nullcone ] = \C[ W ] / ( \momentsmall^{\ast}( x_{ij} ) , \momentsmall^{\ast}( y_{ij} ) ) 
	= \C[ W ] / \C[ W ] \cdot \C[ W ]_+^{\KbigC} \simeq \harmonics( \KbigC ) .
\end{displaymath}
The decomposition of $ \harmonics ( \KbigC ) $ is given in \cite[Th.~3.7.3.1 \& Cor.~3.7.3.6]{Howe.1995.SchurLecture}.  
See also \cite{Gelbart.1974} and \cite{Kashiwara.Vergne.1978}.  
By the same theorem of Kostant, the null cone $ \nullcone $ is normal.  Also, $ \nullconepq $ are normal varieties.
\end{proof}

\subsection{$ U(p, q) \times U(m, n) \  ( m + n \leq \mathoperator{min} ( p, q ) ) $}

Put 
$ \Gbig = U(p, q) $ and $ \Gsmall = U(m, n) $, 
and assume the stable range condition $ m + n \leq p, q $.  We denote:
\begin{gather*}
\begin{array}{ccc}
\Mbig = U(p, q)^2 & \overset{\Delta}{\supset} & \Gbig = U(p, q) \\
\cup & & \cup \\
\Lbig = \Lbigp \times \Lbigq = U(p)^2 \times U(q)^2 & \overset{\Delta}{\supset} & \Kbig = \Kbigp \times \Kbigq = U(p) \times U(q)
\end{array}
\end{gather*}
Here $ \Delta $ denotes the diagonal embedding.  
Similarly, we denote $ \Msmall, \Gsmall, \Lsmall, \Ksmall $ replacing $ p $ and $ q $ by $ m $ and $ n $ respectively.  
Put
\begin{align*}
W 
& = M_{p + q, m + n}(\C) = \left\{ \begin{pmatrix} A & B \\ C & D \end{pmatrix} \mid A \in M_{p, m} , B \in M_{p, n} , 
		C \in M_{q, m} , D \in M_{q, n} \right\} \\
& = \Wp \oplus \Wq = M_{p, m + n}(\C) \oplus M_{q, m + n}(\C) .
\end{align*}
Then $ \KbigC = GL_p(\C) \times GL_q(\C) $ and $ \KsmallC = GL_m(\C) \times GL_n(\C) $ act on $ W $ as 
\begin{align*}
( ( g_1 , g_2 ) , ( h_1 , h_2 ) ) &\cdot \begin{pmatrix} A & B \\ C & D \end{pmatrix} 
= \begin{pmatrix} g_1 A \transpose{h_1} & \transpose{g_1}^{-1} B h_2^{-1} \\ 
		  \transpose{g_2}^{-1} C h_1^{-1} & g_2 D \transpose{h_2} \end{pmatrix} , \\[5pt]
&( g_1 , g_2 ) \in GL_p(\C) \times GL_q(\C) , \quad ( h_1 , h_2 ) \in GL_m(\C) \times GL_n(\C) .
\end{align*}
We fix a Cartan decomposition $ \gbig = \kbig \oplus \sbig $ as 
\begin{align*}
\gbig 
&= \mathfrak{gl}_{p + q}( \C ) = 
	\left( \begin{array}{@{}c|c@{}}
	\rule[-1.3ex]{0pt}{3.3ex} M_p(\C) & 0          \\ \hline
	\rule[-1ex]{0pt}{3.5ex}   0          & M_q(\C) 
	\end{array} \right) 
\oplus 
	\left( \begin{array}{@{}c|c@{}}
	\rule[-1.3ex]{0pt}{3.3ex} 0             & M_{p, q}(\C) \\ \hline
	\rule[-1ex]{0pt}{3.5ex}   M_{q, p}(\C)  & 0            
	\end{array} \right)
= \kbig \oplus \sbig .
\end{align*}
The Cartan decomposition $ \gsmall = \ksmall \oplus \ssmall $ is chosen similarly.  
We identify 
\begin{displaymath}
\sbig = M_{p, q}(\C) \oplus M_{q, p}(\C) , \qquad 
\ssmall = M_{m, n}(\C) \oplus M_{n, m}(\C) .
\end{displaymath}
The moment maps $ \momentbig $ and $ \momentsmall $ are explicitly given by 
\begin{displaymath}
\momentbig : W \ni \begin{pmatrix} A & B \\ C & D \end{pmatrix} \longmapsto ( A \transpose{C} , D \transpose{B} ) \in \sbig , \quad 
\momentsmall : W \ni \begin{pmatrix} A & B \\ C & D \end{pmatrix} \longmapsto ( \transpose{A} B , \transpose{D} C ) \in \ssmall .
\end{displaymath}
These maps are $ \KbigC \times \KsmallC $-equivariant with the trivial $ \KbigC $-action on $ \ssmall $, and 
the trivial $ \KsmallC $-action on $ \sbig $ respectively.  
The maps $ \momentbig $ and $ \momentsmall $ are almost the same.  Therefore we will treat only the map $ \momentsmall $ in the following.  

The map $ \momentsmall $ induces an algebra homomorphism $ \momentsmall : \C[ \ssmall ] \rightarrow \C[ W ] $ by 
\begin{displaymath}
\textstyle
\momentsmall^{\ast}( x_{i, j} ) = \sum_{k = 1}^p a_{k, i} b_{k, j} \, ; \quad 
\momentsmall^{\ast}( y_{i, j} ) = \sum_{l = 1}^q c_{l, j} d_{l, i} , 
\end{displaymath}
where $ ( x_{i, j} ) $ and $ ( y_{i, j} ) $ are coordinate function on $ ( X , Y ) \in M_{m, n}(\C) \times M_{n, m}(\C) = \ssmall $.  
Then the $ GL_p $-invariants (resp. $ GL_q $-invariants) on $ \C[ \Wp ] $ (resp. $ \C[ \Wq ] $) are generated by 
$ \momentsmall( x_{i, j} )^{\ast} $'s (resp.  $ \momentsmall( y_{i, j} )^{\ast} $'s).  
Hence $ \momentsmall : W \rightarrow \ssmall $ is an affine quotient map by $ \KbigC $, 
which is surjective under the condition of the stable range.  Note that $ \momentbig $ is \textit{not} surjective in general.  

Let us consider the null cones
\begin{displaymath}
\nullcone = \momentsmall^{-1}( 0 ) \subset W , \qquad 
\nullconepq = ( \momentsmallpq )^{-1}( 0 ) \subset \Wpq , 
\end{displaymath}
where $ \momentsmallpq : \Wpq \rightarrow \ssmallpq $ is the restriction of $ \momentsmall $ to $ \Wpq $.

\begin{proposition}
\label{prop:UU.nullcone}
Assume the stable range condition $ m + n \leq \min (p, q) $.
\begin{thmenumerate}
\item
\label{item1:prop:UU.nullcone}
The null cone $ \nullcone = \nullconep \times \nullconeq $ is an irreducible normal variety, which is of complete intersection in the set theoretic sense.  
\item
\label{item2:prop:UU.nullcone}
There are finitely many $ \KbigC \times \KsmallC $-orbits in $ \nullcone $, which are completely classified by the ranks of each component in $ \nullconepq $.  
Among them, there is an open dense $ \KbigC \times \KsmallC $-orbit $ \sorbit_{m, n} $, which is a single $ \KbigC $-orbit.    
Moreover, the singular locus is of codimension $ \geq 2 $.
\item
\label{item3:prop:UU.nullcone}
The regular function ring $ \C[ \nullcone ] $ is isomorphic to the harmonics 
$ \harmonics( \KbigC ) \simeq \harmonicsp( \KbigpC ) \otimes \harmonicsq( \KbigqC ) $ as a $ \KbigC \times \KsmallC $-module.  
Its decomposition is explicitly given as 
\begin{align*}
\C[ \nullconep ] &\simeq \harmonicsp( \KbigpC ) 
	\simeq \directsum_{\alpha \in \partition_m , \beta \in \partition_n} 
	\tau_{\alpha \composit \beta}^{(p)}{}^{\ast} \otimes 
		\bigl( \tau_{\alpha}^{(m)}{}^{\ast} \otimes \tau_{\beta}^{(n)} \bigr)  , \\
\C[ \nullconeq ] &\simeq \harmonicsq( \KbigqC ) 
	\simeq \directsum_{\gamma \in \partition_m , \delta \in \partition_n} 
	\tau_{\gamma \composit \delta}^{(q)} \otimes 
		\bigl( \tau_{\gamma}^{(m)} \otimes \tau_{\delta}^{(n)}{}^{\ast} \bigr) , 
\end{align*}
where, for partitions $ \alpha \in \partition_m $ and $ \beta \in \partition_n $, we denote 
\begin{displaymath}
\alpha \composit \beta = \alpha \composit_p \beta 
	= ( \alpha_1 , \alpha_2, \ldots , \alpha_m , 0 , \ldots, 0 , - \beta_n , \ldots , - \beta_1 ) \in \Z^p ; 
\end{displaymath}
and $ \tau_{\lambda}^{(p)} $ is an irreducible finite dimensional representation of 
$ GL_p $ of highest weight $ \lambda $.
\end{thmenumerate}
\end{proposition}

\begin{proof}
Put 
\begin{math}
\sorbitp_{r, s} = \{ ( A , B ) \in M_{p, m + n} \mid \transpose{A} B = 0 , \rank A = r , \rank B = s \} \subset \nullconep .
\end{math}
Then $ \sorbitp_{r, s} $ is a $ \KbigpC \times \KsmallC $-orbit, and 
\begin{math}
\dim \sorbitp_{r, s} = r ( m + p ) + s ( n + p ) + r s - ( r + s )^2 .
\end{math}
It is easy to see that $ \sorbitp_{m, n} \subset \nullconep $ is an open dense orbit and 
\begin{math}
\dim \nullconep = \dim \sorbitp_{m, n} = p ( m + n ) - m n , 
\end{math}
which means $ \nullconep $ is of complete intersection.  The singular locus is given by 
\begin{displaymath}
\closure{\sorbitp_{m - 1, n - 1}} = \coprod_{ r \leq m - 1 , s \leq n - 1 } \sorbitp_{r, s} \subset \nullconep ,
\end{displaymath}
which is of codimension $ \geq 2 $.  The rest of the proof is similar to that of Proposition \ref{prop:OSp.nullcone}.  

The decomposition of harmonics is given in \cite[Th.~2.5.4]{Howe.1995.SchurLecture}.
\end{proof}

\subsection{$ Sp(p, q) \times O^{\ast}(2n) \  ( n \leq \mathoperator{min} ( p, q ) ) $}

Put 
$ \Gbig = Sp(p, q) $ and $ \Gsmall = O^{\ast}(2n) $, 
and assume the stable range condition $ n \leq p, q $.  We denote:
{\small
\begin{gather*}
\begin{array}{ccc}
\Mbig = U(2 p, 2 q) & \supset & \Gbig = Sp(p, q) \\
\cup & & \cup \\
\begin{array}[t]{l@{\,}l}
\Lbig &= \Lbigp \times \Lbigq \\
      &= U(2p) \times U(2q) 
\end{array}
& \supset & 
\begin{array}[t]{l@{\,}l}
\Kbig &= \Kbigp \times \Kbigq \\
      &= Sp(p) \times Sp(q)
\end{array}
\end{array}
\qquad
\begin{array}{ccc}
\Msmall = O^{\ast}(2n)^2 & \overset{\Delta}{\supset} & \Gsmall = O^{\ast}(2n) \\
\cup & & \cup \\
\begin{array}[t]{l@{\,}l}
\Lsmall &= \Ksmall \times \Ksmall \\
        &= U(n)^2 
\end{array}
& \overset{\Delta}{\supset} & \Ksmall = U(n)
\end{array}
\end{gather*}
}
Put
\begin{align*}
W 
& = M_{2 p + 2 q, n}(\C) = \left\{ \begin{pmatrix} A \\ B \end{pmatrix} \mid A \in M_{2 p, n}(\C) , B \in M_{2 q, n}(\C) \right\} \\
& = \Wp \oplus \Wq = M_{2 p, n}(\C) \oplus M_{2 q, n}(\C) .
\end{align*}
The complexified groups $ \KbigC = Sp(2 p, \C) \times Sp(2 q, \C) $ and $ \KsmallC = GL_n(\C) $ act on $ W $ as 
\begin{displaymath}
\begin{pmatrix} A \\ B \end{pmatrix} \mapsto \begin{pmatrix} k_1 A \transpose{g} \\ \transpose{k_2}^{-1} B g^{-1}\end{pmatrix} , \qquad
\begin{array}{r@{\,}l}
( k_1 , k_2 ) &\in Sp(2 p, \C) \times Sp(2 q, \C) \\
g &\in GL_n(\C) 
\end{array}
\end{displaymath}
We fix Cartan decompositions $ \gbig = \kbig \oplus \sbig $ and $ \gsmall = \ksmall \oplus \ssmall $ as 
\begin{align*}
\gbig 
&= \mathfrak{sp}( 2 p + 2 q , \C ) = 
	\left( \begin{array}{@{}c|c@{}}
	\rule[-1.3ex]{0pt}{3.3ex} \mathfrak{sp}(2 p, \C) & 0          \\ \hline
	\rule[-1ex]{0pt}{3.5ex}   0          & \mathfrak{sp}(2 q, \C) 
	\end{array} \right) 
\oplus 
	\left( \begin{array}{@{}c|c@{}}
	\rule[-1.3ex]{0pt}{3.3ex} 0                         & M_{2 p, 2 q} \\ \hline
	\rule[-1ex]{0pt}{3.5ex}   J_q \transpose{M}_{2 p, 2 q} J_p  & 0            
	\end{array} \right) , \\
& \qquad = \kbig \oplus \sbig ,  \quad J_p = \begin{pmatrix} 0 & 1_p \\ - 1_p & 0 \end{pmatrix} \\[10pt]
\gsmall 
&= \mathfrak{o}( 2 n , \C ) = 
	\left( \begin{array}{@{}c|c@{}}
	\rule[-1.3ex]{0pt}{3.3ex} M_n(\C) & 0          \\ \hline
	\rule[-1ex]{0pt}{3.5ex}   0          & - \transpose{M}_n(\C) 
	\end{array} \right) 
\oplus 
	\left( \begin{array}{@{}c|c@{}}
	\rule[-1.3ex]{0pt}{3.3ex} 0           & \Alt_n(\C) \\ \hline
	\rule[-1ex]{0pt}{3.5ex}   \Alt_n(\C)  & 0            
	\end{array} \right)
= \ksmall \oplus \ssmall .
\end{align*}
Thus, we can identify 
$ \sbig = M_{2 p, 2 q}(\C) $ and $ \ssmall = \ssmallp \oplus \ssmallq = \Alt_n(\C) \oplus \Alt_n(\C) $.  
The moment maps are explicitly given by 
\begin{displaymath}
\momentbig : W \ni \begin{pmatrix} A \\ B \end{pmatrix} \longmapsto A \transpose{B} \in \sbig , \qquad 
\momentsmall : W \ni \begin{pmatrix} A \\ B \end{pmatrix} \longmapsto ( \transpose{A} J_p A , \transpose{B} J_q B ) \in \ssmall .
\end{displaymath}
These maps are clearly $ \KbigC \times \KsmallC $-equivariant with the trivial $ \KbigC $-action on $ \ssmall $, 
and the trivial $ \KsmallC $-action on $ \sbig $ respectively.  
They induce algebra morphisms $ \momentbig^{\ast} $ and $ \momentsmall^{\ast} $:
\begin{align*}
\momentbig^{\ast}( z_{i, j} ) &= \textstyle \sum_{k = 1}^n a_{i, k } b_{j, k} , \quad 
	Z = ( z_{i, j} ) \in M_{2 p, 2 q} , \;
	A = ( a_{i, j} ) \in M_{2 p, n} , \, B =( b_{i, j} ) \in M_{2 q, n} , \\
\momentsmall^{\ast}( x_{i, j} ) 
&= \textstyle - \sum_{k = 1}^p a_{k + p, i} a_{k, j} + \sum_{k = 1}^p a_{k, i} a_{k + p, j} , \quad X = ( x_{i, j} ) \in \Alt_n , \\ 
\momentsmall^{\ast}( y_{i, j} ) 
&= \textstyle - \sum_{l = 1}^q b_{l + q, i} b_{l, j} + \sum_{l = 1}^q b_{l, i} b_{l + q, j} , \quad Y = ( y_{i, j} ) \in \Alt_n .
\end{align*}
Classical invariant theory tells us that $ \momentbig^{\ast}( z_{i, j} ) $'s generate the $ GL_n $-invariants $ \C[ W ]^{\KsmallC} $, 
and that $ \momentsmall( x_{i, j} )^{\ast} $'s (resp. $ \momentsmall( y_{i, j} )^{\ast} $'s) generate the 
$ \KbigpC = Sp(2 p, \C) $-invariants $ \C[ \Wp ]^{\KbigpC} $ 
(resp. $ \KbigqC = Sp(2 q, \C) $-invariants $ \C[ \Wq ]^{\KbigqC} $).  
This shows that $ \momentbig $ and $ \momentsmall $ are affine quotient maps.  
By the assumption of the stable range, $ n \leq p, q $, the map 
$ \momentsmall : W \rightarrow \ssmall $ is a surjection, and we have 
\begin{math}
\C[ \ssmall ] \simeq \C[ W ]^{\KbigC} , \;
\C[ \ssmallpq ] \simeq \C[ \Wpq ]^{\KbigpqC} .
\end{math}

Let us define the null cones by 
\begin{displaymath}
\nullcone = \momentsmall^{-1}( 0 ) \subset W , \qquad 
\nullconepq = ( \momentsmallpq )^{-1}( 0 ) \subset \Wpq , 
\end{displaymath}
where $ \momentsmallpq : \Wpq \rightarrow \ssmallpq $ is the restriction of $ \momentsmall $ to $ \Wpq $.

\begin{proposition}
\label{prop:SpOast.nullcone}
Assume the stable range condition $ n \leq p, q $.
\begin{thmenumerate}
\item
\label{item1:prop:SpOast.nullcone}
The null cone $ \nullcone = \nullconep \times \nullconeq $ is an irreducible normal variety, which is of complete intersection 
in the set theoretical sense.
\item
\label{item2:prop:SpOast.nullcone}
There are finitely many $ \KbigC \times \KsmallC $-orbits in $ \nullcone $, which are completely classified by the ranks of each component in $ \nullconepq $.  
The full rank orbit $ \sorbit_{n, n} $ is an open dense $ \KbigC \times \KsmallC $-orbit in $ \nullcone $, which is also a single $ \KbigC $-orbit.      
The singular locus of $ \nullcone $ is of codimension $ \geq 2 $.
\item
\label{item3:prop:SpOast.nullcone}
The regular function ring $ \C[ \nullcone ] $ is isomorphic to the harmonics $ \harmonics( \KbigC ) \simeq \harmonicsp( \KbigpC ) \otimes \harmonicsq( \KbigqC ) $ as a $ \KbigC \times \KsmallC $-module, and we have 
\begin{displaymath}
\C[ \nullconep ] \simeq \harmonicsp( \KbigpC ) 
	\simeq \directsum_{\lambda \in \partition_n} \sigma_{\lambda}^{(p)}{}^{\ast} \otimes \tau_{\lambda}^{(n)}{}^{\ast} ; \quad 
\C[ \nullconeq ] \simeq \harmonicsq( \KbigqC ) 
	\simeq \directsum_{\lambda \in \partition_n} \sigma_{\lambda}^{(q)} \otimes \tau_{\lambda}^{(n)} , 
\end{displaymath}
where $ \sigma_{\lambda}^{(p)} $ denotes an irreducible finite dimensional representation of $ Sp(2 p, \C) $ with highest weight $ \lambda \in \partition_n $.
\end{thmenumerate}
\end{proposition}

\begin{proof}
Put 
\begin{align*}
\sorbitp_{r} &= \{ A \in M_{2 p, n} \mid \transpose{A} J_p A = 0 , \rank A = r \} \subset \nullconep , \quad \\
\sorbitq_{s} &= \{ B \in M_{2 q, n} \mid \transpose{B} J_q B = 0 , \rank B = s \} \subset \nullconeq .
\end{align*}
Then $ \sorbit_{r, s} = \sorbitp_r \times \sorbitq_s $ is a single $ \KbigC \times \KsmallC $-orbit and 
\begin{displaymath}
\dim \sorbitp_r = r ( 2 p + n ) - r^2 - \frac{ r ( r - 1 ) }{2} , \qquad 
\dim \sorbitp_s = s ( 2 q + n ) - s^2 - \frac{ s ( s - 1 ) }{2} .
\end{displaymath}
It is easy to show that $ \sorbitpq_n \subset \nullconepq $ is an open dense orbit and 
\begin{displaymath}
\dim \nullconep = \dim \sorbitp_n = 2 p n  - \frac{n (n - 1)}{2} , \qquad 
\dim \nullconeq = \dim \sorbitq_n = 2 q n  - \frac{n (n - 1)}{2} ,
\end{displaymath}
which tells us that $ \nullcone $ (resp. $ \nullconepq $) is of complete intersection.  
The singular locus is given by 
\begin{displaymath}
\coprod_{ r \leq n - 2} \sorbitpq_r \subset \nullconepq , 
\end{displaymath}
which is of codimension $ \geq 2 $.  
The rest of the proof is similar to that of Proposition \ref{prop:OSp.nullcone}.

The structure of harmonics is given in \cite[Th.~3.8.6.2]{Howe.1995.SchurLecture}.  
\end{proof}

\bigskip


\end{document}